\documentclass[12pt]{article}
\usepackage[intlimits]{amsmath}
\usepackage{amsfonts,amssymb,amscd,amsthm}
\input xy
\input epsf
\xyoption{all}

\newtheorem{theorem}{Theorem}[section]
\newtheorem{lemma}{Lemma}[section]
\newtheorem{proposition}{Proposition}[section]
\newtheorem{remark}{Remark}[section]
\newtheorem{example}{Example}[section]
\newtheorem{definition}{Definition}[section]

\newtheorem{assumption}{Assumption}[section]

\newlength{\rig}
\newlength{\rigg}
\newlength{\hei}
\newcommand{\dokaend}{\hfill$\square$ \vskip6truept}

\newcommand{\fgr}[3]{
\setlength{\hei}{7.8cm} \setlength{\rig}{0.1\textwidth}
\setlength{\rigg}{0.8\textwidth}
\begin{figure}
\rule{\rig}{0in}
\epsfxsize=13cm
\epsffile[1 1 953 638]{#1} \caption{#3}\label{#2}
\end{figure}
}

\numberwithin{equation}{section}
\newcommand{\func}[1]{{\rm #1} \,}

\newcommand{\limfunc}[1]{{\rm{#1}}}
\def\stackunder#1#2{\mathrel{\mathop{#2}\limits_{#1}}}%

\title{\bf Index Defects in the Theory of Non-local
Boundary Value Problems and the $\eta$-Invariant}

\author{{\bf A.~Yu.~Savin and B.~Yu.~Sternin}}

\begin{document}
\maketitle

\begin{abstract}
The paper deals with elliptic theory of boundary value problems on manifolds
whose boundary is represented as a covering space. We compute the index for a
class of non-local boundary value problems on such manifolds. For a non-trivial
covering, the index defect of the Atiyah--Patodi--Singer boundary value problem
is computed. Poincar\'{e} duality in $K$-theory of the corresponding manifolds
with singularities is obtained.
\end{abstract}

\section*{Introduction}
\addcontentsline{toc}{section}{Introduction}

This paper deals with boundary value problems for elliptic operators on a
manifold whose boundary is the total space of a finite-sheeted covering. On
such manifolds, we consider boundary value problems for operators that do not
satisfy the Atiyah--Bott condition (i.e., have no well-posed classical boundary
value problems). Recall that this condition does not hold, in particular, for
the  Hirzebruch and Dirac operators as well as some other related geometric
operators.

We consider the following two classes of boundary value problems.

\textbf{1. Non-local boundary value problems}. Let $M$ be a smooth manifold
such that the boundary $\partial M$  is a finite-sheeted covering with
projection $ \pi :\partial M\longrightarrow X.$ Then there is an isomorphism
$$ C^\infty \left( \partial M\right)
\stackrel{\beta }{\simeq }C^\infty \left( X,\pi _{!}1\right)
$$
between the space $C^\infty \left( \partial M\right) $ of smooth functions on
$\partial M$ and the space of sections of the vector bundle $\pi _{!}1\in
\limfunc{Vect}\left( X\right) $ on the base of the covering. Here $\pi _{!}1$
is the direct image of the trivial line bundle.

For a scalar elliptic operator $D$ on $M$, the simplest non-local boundary
value problem of the type considered in this paper is
\begin{equation}
\left\{
\begin{array}{l}
Du=f, \\
B\beta \left. u\right| _{\partial M}=g.
\end{array}
\right.  \label{sisa2}
\end{equation}
Here $u$ and $f$ are functions on $M$, $g$ is a function on $X$, and the
operator $B$ of boundary conditions acts also on $X.$ In terms of the original
manifold $M$, the boundary conditions in (\ref{sisa2}) are {\em non-local},
since they relate the values of functions at distinct points of $M$.

We prove a finiteness theorem and in the case of regular coverings obtain an
index formula for this class of non-local boundary value problems. Without
going into detail at the moment, let us mention two essential features of the
theory.

First, in the proof of the index theorem we embed our manifolds in  the
classifying space of a finite group, while in the classical index theorem  it
is suffices to use embeddings in $\Bbb{R}^N$.

Second, the analogue of the Atiyah--Singer difference element for a non-local
boundary value problem is an element of the $K$-group of a {\em non-commutative
$C^{*}$-algebra} associated with the cotangent bundle and the covering. Recall
that in the classical index theorem it suffices to use topological $K$-theory.

The index formula of this paper is given in a form resembling the $K$-theoretic
statement of the Atiyah--Singer theorem. Local index formulae will appear
elsewhere.\vspace{1mm}

\textbf{2. Spectral problems on manifolds with a covering}. The first
generalization of classical boundary value problems that is free of the
Atiyah--Bott obstruction is due to Atiyah, Patodi, and Singer \cite{APS1}. For
a class of first-order elliptic operators, one has so-called \emph{spectral
boundary value problems} denoted by $\left(D,\Pi_{+}\right).$ Spectral boundary
value problems enjoy the Fredholm property. However, their index is not
determined by the principal symbol of $D$.

Interesting invariants arise if the boundary has the structure of a covering.
Here we consider a class of elliptic operators that are lifted from the base of
the covering in a neighbourhood of the boundary. In this case, the principal
symbol of an elliptic operator $D$ defines an element
\[
\left[ \sigma \left( D\right) \right] \in K^0\left( \overline{T^{*}M}^\pi
\right)
\]
in the $K$-group of the singular space $\overline{T^{*}M}^\pi $ obtained from
the cotangent bundle $T^{*}M$ if we identify all points in each fiber of the
covering (for details, see Section \ref{defe3}). The element $\left[ \sigma
\left( D\right) \right] $ has a topological index
\[
\limfunc{ind}_t\left[ \sigma \left( D\right) \right] \in \Bbb{Q}/n%
\Bbb{Z},
\]
where $n$ is the number of sheets. However, the analytical and the topological
index coincide only for trivial coverings. For a general covering, we obtain
the \emph{index defect formula}
\begin{equation}
\func{mod}n\text{-}\limfunc{ind}\left( D,\Pi _{+}\right) -\limfunc{ind}%
_t\left[ \sigma \left( D\right) \right] =\eta \left( \left.
D\right| _{X}\otimes 1_{n-\pi _{!}1}\right) \in \Bbb{Q}/n\Bbb{Z}.
\label{def2}
\end{equation}
The index defect (the difference between the analytical index modulo $n$ and
the topological index) is equal to the relative Atiyah--Patodi--Singer
$\eta$-invariant of the restriction of $D$ to the boundary with coefficients in
the flat bundle $\pi _{!}1.$ For a trivial covering, the relative $\eta
$-invariant is zero, and the index defect formula becomes the index formula
\[
\func{mod}n\text{-}\limfunc{ind}\left( D,\Pi _{+}\right) =\limfunc{ind}%
_t\left[ \sigma \left( D\right) \right]
\]
due to Melrose and Freed \cite{FrMe1} (see also \cite{Hig1,Zha1,Bot1,Ros2}).
However, our proof is new even in this case. It is interesting to note that the
main step in the proof is to realize the fractional analytic invariant
\begin{equation}
\func{mod}n\text{-}\limfunc{ind}\left( D,\Pi _{+}\right) -\eta \left( \left.
D\right| _{X}\otimes 1_{n-\pi _{!}1}\right) \in \Bbb{Q}/n\Bbb{Z}
\label{in7}
\end{equation}
as the index of some non-local boundary problem of the form (\ref{sisa2}) (in a
suitable elliptic theory with coefficients).

There is also a deeper relation between the two elliptic theories described in
Subsec.~1 and 2. \vspace{1mm}

{\bf 3. Poincar\'{e} isomorphism and duality.} We establish {\em Poincar\'{e} isomorphisms} on
the singular spaces $\overline{T^*M}^\pi$ and $\overline{M}^\pi$. For the identity covering
$\pi=Id,X=\partial M$, these isomorphisms  are just the well-known isomorphisms (e.g., see
\cite{Con1,Kas3,MePi2})
\begin{equation}
\label{alox2}
 K^0(T^*M)\simeq K_0(M,\partial M), \quad
K^0(T^*(M\setminus\partial M))\simeq K_0(M).
\end{equation}
(For non-compact spaces, we use $K$-theory with compact supports.) In contrast
to the smooth case, the Poincar\'{e} isomorphisms for singular spaces relate
the $K$-groups of a commutative algebra of functions to those of a dual
non-commutative algebra. They are defined on the elements as {\em
quantizations}, i.e., take symbols to  operators. More precisely, the analogue
of the first isomorphism in (\ref{alox2}) is defined in terms of the operators
described in Subsec.~2, while in the second case one uses non-local problems
introduced in Subsec.~1.
%It should be mentioned that in the $C^*$-algebra $K$-theory one knows an abstract form of
%duality (see \cite{Hig2,Pash1}).

Let us outline the contents of the paper. The first section contains the
definition of the class of non-local boundary value problems on manifolds with
a covering on the boundary and a proof of the Fredholm property. The index
formula is obtained in Sec.~2. By way of example, we define a non-local
boundary value problem for the Hirzebruch operator on a manifold with
reflecting boundary. In Sec.~\ref{fama1}, we give  the homotopy classification
of non-local problems.  The index defect formula  (\ref{def2}) is proved in
Sec.~\ref{defe3}. This is one of the central results of the paper. Section
\ref{exa4} contains applications to the computation of the fractional part of
the $\eta $-invariant. It is also shown that the invariant (\ref{in7}) can be
computed by the Lefschetz formula. Poincar\'{e} isomorphisms in $K$-theory of
the singular spaces corresponding to manifolds whose boundary bears the
structure of a covering are constructed in the last two sections.

There are other interesting classes of non-local boundary value problems
arising if the projection has singularities (e.g., the projection on the
quotient by a non-free action of a finite group). Index theory of such boundary
value problems is apparently related to index theory on orbifolds (\cite{Kaw1,
Far1}). Our approach is advantageous in that if the base of the covering is
smooth, then there are no additional analytic and topological difficulties
related to the singularities of the covering. More general classes of non-local
boundary value problems (e.g., see \cite{Ant1}) are beyond the scope of this
paper.\vspace{1mm}

\textbf{Acknowledgements}. The results were announced at the conferences
``Spring School 2001'' in Potsdam, Germany, ``Topology, analysis, and related
topics'' in Moscow, 2001, and at the International Congress of Mathematicians
in Beijing, 2002.  We are grateful to V.E.~Nazaikinskii and V.~Nistor for
helpful discussions. The work was supported in part by RFBR grants Nos.
03-02-16336, 02-01-00118, and 02-01-00928. We are grateful to A.S.~Mishchenko,
for finding an error in the original version of the paper and  to the referee
for helpful remarks.

\section{Non-local boundary value problems}\label{para1}

\textbf{1}. \textbf{Coverings and non-local operators}. Let $Y$ be a finite
covering over a manifold $X$ with projection $\pi :Y\longrightarrow X.$ The
projection defines the direct image mapping
\[
\pi _{!}:\limfunc{Vect}\left( Y\right) \longrightarrow \limfunc{Vect}\left(
X\right)
\]
that takes each vector bundle $E\in \limfunc{Vect}\left(Y\right) $ to the
bundle
\[
\text{ }\pi _{!}E\in \limfunc{Vect}\left( X\right) ,\quad \left( \pi
_{!}E\right) _x=C^\infty \left( \pi ^{-1}\left( x\right) ,E\right) ,\quad x\in
X.
\]
This clearly gives an isomorphism
$\beta _E:C^\infty \left( Y,E\right) \stackrel{\simeq }{\longrightarrow }%
C^\infty \left( X,\pi _{!}E\right)$ of section spaces on $Y$ and $X$, while
permits one to identify operators defined on the total space and on the base.
More precisely, the direct image
\[
\pi _{!}D=\beta _ED\beta _E^{-1}:C^\infty \left( X,\pi _{!}E\right)
\longrightarrow C^\infty \left( X,\pi _{!}E\right)
\]
of a differential operator
\[
D:C^\infty \left( Y,E\right) \longrightarrow C^\infty \left( Y,E\right)
\]
on $Y$ is a differential operator. However, the following example shows that
the inverse image
\begin{equation}
\pi ^{!}D^{\prime }=\beta _E^{-1}D^{\prime }\beta _E:C^\infty \left(
Y,E\right) \longrightarrow C^\infty \left( Y,E\right)  \label{nona}
\end{equation}
of a differential operator $D^{\prime } $ on $X$ may well be a \emph{non-local
operator}. (It is not even pseudolocal.)

\begin{example}
\emph{For the trivial covering}
\[
Y=\stackunder{n\;\;{\emph{copies}}}{\underbrace{X\sqcup X\sqcup ...\sqcup
X}}\longrightarrow X
\]
\emph{\ and the trivial bundle }$E=\Bbb{C},$\emph{\ we have }$\pi
_{!}E=\Bbb{C}^n.$\emph{\ The direct image}
\[
\pi _{!}D=\limfunc{diag}\left( \left. D\right| _{X_1},...,\left. D\right|
_{X_n}\right) :C^\infty \left( X,\Bbb{C}^n\right) \longrightarrow C^\infty
\left( X,\Bbb{C}^n\right)
\]
{\em of a differential operator on }$Y$\emph{\ is always a diagonal operator,
and hence the inverse image of a non-diagonal operator can not be a
differential operator. The off-diagonal entries produce non-local operators on}
$Y$\emph{, since they interchange the values of functions on different leaves
of the covering.}
\end{example}

\textbf{2. Non-local boundary value problems}. Let $M$ be a smooth compact
manifold with boundary $\partial M.$ Suppose that the boundary is a covering
space over a smooth closed manifold $X$ with projection
\[
\pi :\partial M\longrightarrow X.
\]
We fix a collar neighbourhood $\partial M\times \left[0,1\right) $ of the
boundary. The normal coordinate will be denoted by $t$.

For a smooth function $u\in C^\infty \left( M\right) $,  let
\[
j_{\partial M}^{m-1}u=\left( \left. u\right| _{\partial M},\left. -i\frac
\partial {\partial t}u\right| _{\partial M},\ldots ,\left. \left( -i\frac
\partial {\partial t}\right) ^{m-1}u\right| _{\partial M}\right)
\]
be the restriction of its $(m-1)$st jet in the normal direction to the
boundary. The operator $j_{\partial M}^{m-1}$ is continuous in the Sobolev
spaces
\[
j_{\partial M}^{m-1}:H^s\left( M\right) \longrightarrow
\bigoplus_{k=0}^{m-1}H^{s-1/2-k}\left( \partial M\right) ,\qquad s>m-1/2.
\]
Throughout the paper we assume that for vector bundles $E$ on manifolds with
boundary there are given isomorphisms $p^*(\left.E\right|_{\partial M})\simeq
\left. E\right|_{\partial M\times[0,1]}$ in the collar neighbourhood of the
boundary, where $p:\partial M\times[0,1]\to
\partial M$ is the natural projection. In this case, the normal jet of a section of
$E$ is also well defined.
\begin{definition}
A non-local boundary value problem\emph{\ for a differential operator}
\[
D:C^\infty \left( M,E\right) \longrightarrow C^\infty \left( M,F\right)
\]
\emph{of order }$m$\emph{\ is a system of equations}
\begin{equation}
\left\{
\begin{array}{ll}
Du=f, & u\in H^s\left( M,E\right) ,f\in H^{s-m}\left( M,F\right) ,%
\vspace{2mm} \\
B\beta _Ej_{\partial M}^{m-1}u=g,\qquad & g\in H^\delta \left( X,G\right) ,
\end{array}
\right.  \label{non1}
\end{equation}
\emph{\ where the boundary condition is defined by a pseudodifferential operator}
\[
B:\bigoplus_{k=0}^{m-1}H^{s-1/2-k}\left( X,\pi _{!}\left. E\right|
_{\partial M}\right) \longrightarrow H^\delta \left( X,G\right)
\]
\emph{\ on }$X.$\emph{\  We assume that the component}
\[
B_k:H^{s-1/2-k}\left( X,\pi _{!}\left. E\right| _{\partial M}\right)
\longrightarrow H^\delta \left( X,G\right)
\]
\emph{of }$B$\emph{\ has the order} $s-1/2-k-\delta. $
\end{definition}

\begin{remark}
{\em One can also consider problems similar to (\ref{non1}) in which the
components of the vector function $g$ belong to Sobolev spaces of different
orders. The case in which all components have the same order $\delta$ is more
convenient and can always be achieved by order reduction. }
\end{remark}

Note that for the identity covering $\pi =Id,$ $X=\partial M$, problem
(\ref{non1}) is just a classical boundary value problem
(e.g., see \cite{Hor3}).\vspace{2mm}%

\textbf{3. Relation to classical boundary value problems. Finiteness theorem}.
Note that $\partial M\times \left[ 0,1\right) $  is also a covering with
projection
\[
\pi \times 1:\partial M\times \left[ 0,1\right) \longrightarrow X\times
\left[ 0,1\right) .
\]
The induced isomorphism of function spaces will be denoted by
\[
\beta _E^{\prime }:C^\infty \left( \partial M\times \left[ 0,1\right)
,E\right) \longrightarrow C^\infty \left( X\times \left[ 0,1\right) ,\pi
_{!}E\right) .
\]
The non-local problem $\left( D,B\right) $ can be represented in a
neighbourhood of the boundary as the inverse image of the classical boundary
value problem
\[ \left(
\begin{array}{cc}
\beta _F^{\prime } & 0 \\
0 & 1
\end{array}
\right) \circ \left( %
\begin{array}{c}
D \\
B\beta _Ej_{\partial M}^{m-1}
\end{array}
\right) \circ \left( \beta _E^{\prime }\right) ^{-1}=\left( %
\begin{array}{c}
\beta _F^{\prime }D\left( \beta _E^{\prime }\right) ^{-1}\vspace{2mm} \\
Bj_X^{m-1}
\end{array}
\right).
\]
More specifically, this is the boundary value problem
\begin{equation}
\left(
\begin{array}{c}
\beta _F^{\prime }D\left( \beta _E^{\prime }\right) ^{-1}\vspace{2mm} \\
Bj_X^{m-1}
\end{array}
\right) :C^\infty \left( X\times \left[ 0,1\right) ,\pi _{!}E\right)
\longrightarrow
\begin{array}{c}
C^\infty \left( X\times \left[ 0,1\right) ,\pi _{!}F\right) \\
\oplus \\
C^\infty \left( X,G\right) ,
\end{array}
\label{cla2}
\end{equation}
for the differential operator $\left( \pi \times 1\right) _{!}D=\beta
_F^{\prime }D\left( \beta _E^{\prime }\right) ^{-1}$ on the cylinder $X\times
\left[ 0,1\right) .$ In local coordinates, this operator is represented by a
diagonal matrix with elements acting on different leaves of the covering.
Problem (\ref{cla2}) will be denoted by $\left( \left( \pi \times 1\right)
_{!}D,B\right)$ for short.

We point out that the  \emph{classical boundary value problem}(\ref{cla2})%
\emph{\ is defined only in a neighbourhood of the boundary}, since the covering
is defined only near the boundary.

\begin{definition}
\emph{Problem }$\left( D,B\right) $\emph{\ is said to be }elliptic\emph{ if
}$D$\emph{\ is elliptic and} $\left( \left( \pi \times 1\right) _{!}D,B\right)
$\emph{\ is elliptic, i.e.,  satisfies the Shapiro--Lopatinskii
condition (e.g., see \cite{Hor3}).}%
\end{definition}

The proof of the following finiteness theorem is standard.

\begin{theorem}
\label{fini}An elliptic boundary value problem $\mathcal{D}=\left(
D,B\right) $ defines a Fredholm operator.
\end{theorem}

\emph{Proof}. Let  $D^{-1}$ be the parametrix of $D$ in the interior of the
manifold. Similarly, the parametrix of the classical boundary value problem on
$X\times \left[0,1\right) $ will be denoted by $(L,K).$ They can be pasted
together globally on $M$ by the formula
\[
\mathcal{D}^{-1}=\left(\psi _1D^{-1}\varphi _1+\psi _2\left( \pi \times 1\right) ^{!}L\varphi
_2,\psi_2 K\right).
\]
Here
\begin{equation}
\begin{array}{c}
\varphi_1+\varphi_2=1,\quad \psi_j\varphi_j=\varphi_j\\
\psi_1=0 \text{ near the boundary }, \psi_2=0 \text{ far from the boundary.}
\end{array}
\label{usa2}
\end{equation}
Furthermore, $\psi_2$ is assumed to be constant in the fiber of $\pi\times 1$.
Obviously, $\mathcal{D}^{-1}$ is a two-sided parametrix of $\mathcal{D}$. The
proof is complete. \dokaend

\section{The index of non-local problems\label{hom2}
}

In the previous section, non-local boundary value problems were represented
near the boundary in terms of equivalent classical boundary value problems.
Therefore, we can apply well-known topological methods (e.g., see
\cite{SaScS4}) to compute the index of non-local boundary value
problems.\vspace{1mm}

\textbf{1. Reduction to zero-order operators}. We introduce a class of
operators that are non-local in a neighbourhood of the boundary. A linear
operator
\[
D:C^\infty \left( M,E\right) \longrightarrow C^\infty \left( M,F\right)
\]
will be called an \emph{admissible operator of order $m$} if it can be
represented modulo operators with smooth kernels as
\begin{equation}
D=\psi _1D^{\prime }\varphi _1+\psi _2\left( \pi \times 1\right) ^{!}D^{\prime
\prime }\varphi _2  \label{dek1}
\end{equation}
for cutoff functions $\varphi _{1,2},\psi _{1,2}$ as in the proof of Theorem
\ref{fini}, a pseudodifferential operator $D^{\prime }:C^\infty \left(
M,E\right) \rightarrow C^\infty \left( M,F\right),$ and an operator
\[
D^{\prime \prime }:C^\infty \left( X\times \left[ 0,1\right) ,\pi _{!}E\right)
\longrightarrow C^\infty \left( X\times \left[ 0,1\right) ,\pi _{!}F\right)
\]
that is a sum of a pseudodifferential operator with compactly supported kernel
on $X\times(0,1)$ and a differential operator
\begin{equation}
\sum_{k=0}^mD_k\left( t\right) \left( -i\frac \partial {\partial t}\right)
^{m-k}  \label{alef1}
\end{equation}
with respect to the normal variable $t$.

Here the $D_k\left( t\right) $\ are smooth families of pseudodifferential
operators on $X$\ of order $k$ and $D_0\left( t\right) $ is induced by a vector
bundle isomorphism.

To this class of operators, one can extend the notion of ellipticity, the
statement of non-local boundary value problems, and the finiteness theorem (cf.
a similar generalization in \cite{Hor3} for the classical case). In particular,
the \emph{symbol of an admissible operator} is a pair $(\sigma_M,\sigma_X),$
where $\sigma_M:p^*E\to p^*F$ is defined over $M\setminus(\partial
M\times[0,\varepsilon))$ ($p:S^*M\to M$ is the natural projection) and
$\sigma_X:p^*_0(\pi_!E)\to p^*_0(\pi_!F)$ ($p_0:S^*(X\times[0,1])\to
X\times[0,1]$) is defined over $X\times[0,1]$. Moreover, the symbols are smooth
and satisfy the compatibility condition
$$
({\pi_0})_!\left.\sigma_M\right|_{\partial M\times(\varepsilon,1)}= \sigma_X,
$$
where the direct image is induced by the natural projection $\pi_0:\left.T^*M\right|_{\partial
M\times(\varepsilon,1)}\to T^*(X\times (\varepsilon,1)).$

\begin{example}
\label{exo1}\emph{Let }$E\in \limfunc{Vect}\left( M\right) $\emph{\ be a vector
bundle. Suppose that  its direct image over $U_{\partial M}$ is decomposed as a
sum of two subbundles}
\[
\pi _{!}\left. E\right| _{U_{\partial M}}=E_{+}\oplus E_{-},\qquad E_{\pm }\in
\limfunc{Vect}\left( X\times[0,\varepsilon)\right) .
\]
\emph{Consider the operator} $D_\pm:C^\infty(M,E)\to C^\infty(M,E)$ {\em given
by}
\begin{equation}
D_{\pm }=\psi _2\left( \pi \times 1\right) ^{!}\left[ \left( -i\frac
\partial {\partial t}+i{\Lambda _{X,E_{+}}}\right) \oplus \left( i\frac
\partial {\partial t}+i{\Lambda _{X,E_{-}}}\right) \right] \varphi
_2+\psi _1i{\Lambda _M}\varphi _1.  \label{kru2}
\end{equation}
\emph{(Here }$\Lambda $%
\emph{\ stands for first-order pseudodifferential operators with principal
symbol $|\xi|$ on the corresponding manifolds, and the cutoff functions are
chosen as before.) This formula defines an admissible elliptic operator. We
equip it with the Dirichlet boundary condition}
\[
P_{E_{-}}\beta_E \left. u\right| _{\partial M}=g\in C^\infty \left( X,E_{-}\right) ,
\]
\emph{where }$P_{E_{-}}:\pi _{!}\left. E\right| _{\partial
M}\rightarrow \pi _{!}\left. E\right|
_{\partial M}$\emph{\ is a projection onto the subbundle }$%
E_{-}.$\emph{\ Denote this boundary value problem } \emph{\ by }
$\mathcal{D}%
_{\pm }.$\emph{\ By analogy with the classical case (e.g., see \cite{Hor3}),
one proves that the index of this boundary value problem is zero}.

\emph{For example, let} $E_+= \pi _{!}\left. E\right| _{\partial M}$ and
$E_{-}=0.$ \emph{Then the operator (\ref{kru2}), which will be denoted by
}$D_{+}$\emph{, is Fredholm without any boundary condition.}
\end{example}

\begin{remark}\emph{
Just as in the classical elliptic theory on a closed manifold
(see~\cite{AtSi1}), there are two equivalent definitions of homotopy of
non-local elliptic problems. First, one can say that two problems are homotopic
if they can be connected by a family of non-local elliptic problems  continuous
in the operator norm (in some given pair of Sobolev spaces). Second, two
problems are said to be homotopic if there exists a continuous homotopy of
their principal symbols (preserving ellipticity). The equivalence of the two
definitions is based  on the smoothing of continuous homotopies and the
standard norm estimates modulo compact operators, e.g., see \cite{KoNi1}.}
\end{remark}

Let\footnote{Later on, by $\limfunc{Ell}^m\left( M,\pi \right)$ we denote also
the corresponding Grothendieck groups for closed manifolds and for manifolds
with boundary with projection $\pi$ defined, possibly, on an open subset. Which
group is meant is always clear from the context.} $\limfunc{Ell}^m\left( M,\pi
\right), $ $m\geq 1$,  be the Grothendieck group of the semigroup of homotopy
classes of elliptic boundary value problems for admissible operators of order
$m$ modulo boundary value problems of the form $\mathcal{D}_{\pm }\circ
D_{+}^{m-1}.$

The group of stable homotopy classes of zero-order admissible elliptic
operators is denoted by $\limfunc{Ell}^0\left( M,\pi \right) $. Recall that
stabilization is taken modulo trivial operators. In this case, by  trivial
operators we mean operators induced by vector bundle isomorphisms. It should be
noted that elliptic operators of order zero do not require boundary conditions,
since near the boundary they are induced by vector bundle isomorphisms.

Just as in the classical theory (see \cite{SaScS4} or \cite{Hor3}), the order
of a non-local boundary value problem can be reduced to zero by stable
homotopies. More precisely, the following theorem holds.

\begin{theorem}[\rm order reduction]\label{th0q}
 The composition with the operator $D_{+}$ \emph{(}with
coefficients in vector bundles\emph{)} induces an isomorphism
\[
\begin{array}{ccc}
\times D_{+}^m:\limfunc{Ell}^0\left( M,\pi \right) & \longrightarrow
& \limfunc{Ell}^m\left( M,\pi \right) ,\vspace{2mm} \\
\qquad\left[ D\right] & \mapsto & \left[ D\circ D_{+}^m\right] .%
\end{array}
\]
\end{theorem}

The \emph{proof} of this result is a straightforward generalization of the
corresponding proof in the classical case (see \cite{SaScS4}) and hence is
omitted. \dokaend

\begin{remark}
{\em \label{remk} Let us explicitly describe order reduction, i.e., the mapping
$\left(\times D^m_+\right)^{-1}$, in the important special case of boundary
value problems
\[
\left\{
\begin{array}{lc}
Du=f, &  \\
P\beta _E\left( \left. u\right| _{\partial M}\right) =g, & \quad g\in C^\infty(X,{\rm Im} P),
\end{array}
\right.
\]
for a first-order admissible  operator $D:C^\infty \left( M,E\right)
\rightarrow C^\infty \left( M,F\right) $ that admits a decomposition
\[
\left( \pi \times 1\right) _{!}\left( \left. D\right| _{U_{\partial M}}\right)
=\Gamma \left( \frac \partial {\partial t}+A\left( t\right) \right)
\]
in a neighbourhood of the boundary, where $A\left( t\right) $ is a smooth
operator family on $X$ and $\Gamma:\pi_!E|_{\partial M}\to \pi_!F|_{\partial
M}$ is a vector bundle isomorphism. The boundary condition is defined by the
projection
 $P$\ in the bundle $\pi _{!}\left( \left.E\right|
_{\partial M}\right).$ We assume for simplicity that the symbol $a(x,\xi)$ of
$A\left( 0\right) $ is symmetric and additionally satisfies $a^{*}a=\left| \xi
\right| ^2.$ We assume that  $P$ is also  symmetric.

Let $L_{+}\left( A\left( 0\right) \right)\in {\rm Vect}(S^*X) $ be the
Calder\'on bundle. For our first-order operator,  this is the bundle over
$S^{*}X$ generated by eigenvectors of $a\left(x, \xi \right) $ with positive
eigenvalues.

The ellipticity condition for $\left( D,P\right) $ requires that $P$ define an
isomorphism
\[
L_{+}\left( A\left( 0\right) \right) \stackrel{P}{\longrightarrow }p _0^{*}%
\func{Im}P,\qquad p _0:S^{*}X\rightarrow X,
\]
of subbundles. Consider the principal symbol of our operator on the boundary:
\[
\sigma \left( \frac \partial {\partial t}+A\left( 0\right) \right) =i\tau +a\left(x, \xi \right)
\]
(here $\tau$ is dual to $t$). The  linear homotopy
\[
\left( 1-\varepsilon \right) \left( i\tau +a\left(x, \xi \right) \right)
+\varepsilon \left( 2P(x)-1\right) ,\quad \varepsilon \in \left[ 0,1\right] ,
\]
is a homotopy of elliptic symbols for $\tau ^2+\xi ^2=1$ provided that the
ellipticity condition for $(D,P)$ is satisfied. Furthermore, at the end of the
homotopy (for $\varepsilon=1$) the symbol does not depend on the cotangent
variables. Let us treat the homotopy of elliptic symbols on $X$ as an elliptic
symbol on $X\times \left[ 0,1\right] $. Then the symbol of $D$ and the homotopy
taken together define the symbol of an admissible elliptic operator of order
zero on the manifold $M$\ with  $\left[ 0,1\right] \times \partial M$ attached.

This zero-order symbol can be transferred to $M$ by an obvious diffeomorphism
$M\simeq M\cup_{\partial M}\left(\left[0,1\right]\times\partial M\right)$ that
is equal to identity far from the boundary. One can show (cf. \cite{SaScS4})
that  the element defined by this symbol (operator) is precisely the image of
the problem
 $\left( D,P\right) $ under the order reduction mapping  $%
\left( \times D_{+}\right) ^{-1}$ of Theorem \ref{th0q}. }
\end{remark}

\textbf{2. Admissible operators on closed manifolds}. Let $\overline{U}$ be a
codimension zero submanifold of some closed manifold $M$. We assume that
$\overline{U}$ is a covering space
\[
\overline{U}\stackrel{\pi }{\longrightarrow }\overline{Y}
\]
with smooth base $\overline{Y}.$ Let $U$ and $Y$ be the corresponding sets of
interior points (we allow $\overline{U}$ to have a boundary). Then {\em scalar
admissible operators} on  $M$ are by definition operators of the form
\[
D=D^{\prime }+\psi \left( \pi ^{!}D^{\prime \prime }\right) \varphi ,
\]
where $D^{\prime }$ is a pseudodifferential operator on $M,$  $D^{\prime \prime }$ is a
pseudodifferential operator on $Y$ acting on sections of  $\pi _{!}1\in \limfunc{%
Vect}\left( Y\right) ,$ and the cutoff functions $\varphi$ and $\psi $ are
supported in $U.$

In the non-scalar case, we consider operators acting in the spaces slightly
more general than section spaces of vector bundles.

Namely, consider triples $\left( E,E_0,\alpha \right)$ defined by vector bundles
$$E\in \limfunc{Vect}\left( V\right) ,\;\; E_0\in
\limfunc{Vect}\left( Y\right)
$$
(here we fix a neighbourhood $V\subset M$ of $M\backslash U$ such that if a
point lies in $U\cap V$ then the entire fiber containing this point also lies
in $U\cap V$) and a vector bundle isomorphism
$$
\pi _{!}\left. E\right|_{U\cap V}\stackrel{\alpha}\simeq E_0|_{\pi\left(U\cap V\right)} $$ on
$\pi\left(U\cap V\right)$.

Let ${\rm Vect}(M,\pi)$ be the set of isomorphism classes of such triples. Here
two triples $(E,E_0,\alpha),(F,F_0,\gamma)$ are \emph{isomorphic} if the vector
bundles are pairwise isomorphic,
$E\stackrel{a}{\simeq}F,E_0\stackrel{b}{\simeq}F_0$, and the isomorphisms are
compatible: $\gamma(\pi_!a)=b\alpha$.

The  linear \emph{space of sections} corresponding to the triple
 ${\cal E}=\left(E,E_0,\alpha\right)$ is defined as
\[
C^\infty \left( M,{\cal E}\right) =\left\{ \left( u,v\right) \;\left|
\begin{array}{c}
u\in C^\infty \left( V,E\right) ,v\in C^\infty \left( Y,E_0\right) , \\
\alpha \beta _E\left( \left. u\right| _{U\cap V}\right) =\left. v\right|%
_{\pi \left( U\cap V\right) }%
\end{array}
\right. \right\} \subset C^\infty \left( V,E\right) \oplus C^\infty \left( Y,E_0\right) .
\]
For the identity covering, ${\cal E} $ defines a vector bundle on $M$ obtained
by clutching $E$ with $E_0$  by the transition function $\alpha$, and $C^\infty
\left( M,{\cal E}\right) $ is just the space of sections of $\mathcal{E}$.

The space $C^\infty \left( M,{\cal E}\right)$ is generated by the subspaces
\[
C_0^\infty \left( V,E\right) ,C_0^\infty \left( Y,E_0\right) \subset C^\infty \left( M,{\cal
E}\right)
\]
of compactly supported sections. More precisely, the first embedding takes $u$
to the pair $(u,\widetilde{\beta_E \left. u\right| _{U\cap V}})$, where the
tilde stands for the  extension of a function by zero at the points where the
function was not originally defined. Similarly, the second embedding takes  $v$
to the pair $(\widetilde{\beta _E^{-1}\left. v\right| _{\pi(U\cap V)}},v).$ Now
non-local operators acting in spaces $C^\infty \left( M,{\cal E}\right) $  can
readily be defined by analogy with the scalar case. Namely, an \emph{admissible
operator of order $m$} is an operator
\[
D:C^\infty \left( M,{\cal E}\right) \longrightarrow C^\infty \left( M,{\cal
F}\right)
\]
that is equal, modulo operators with smooth kernel, to
\begin{equation}
D=D_1\varphi _1+D_2\varphi _2,  \label{omo1}
\end{equation}
where
\[
D_1:C_0^\infty \left( V,E\right) \rightarrow C_0^\infty \left( V,F\right) ,\quad D_2:C_0^\infty
\left( Y,E_0\right) \rightarrow C_0^\infty \left( Y,F_0\right)
\]
are $m$th-order pseudodifferential operators with compactly supported kernels.
Here we assume that the cutoff function $\varphi _1$ is zero in some
neighbourhood of $M\backslash V$ and $\varphi _2$ is zero in a neighbourhood of
$M\backslash U.$

The \emph{symbol of an admissible operator} is a pair
$\left(\sigma_M,\sigma_Y\right)$ of usual elliptic symbols
\[
\sigma _M:p_M^{*}\left. E\right| _{M\backslash U}\longrightarrow p_M^{*}\left. F\right|
_{M\backslash U},\quad \sigma _Y:p_Y^{*}\left. E_0\right| _{\overline{Y}}\longrightarrow
p_Y^{*}\left. F_0\right| _{\overline{Y}},
\]
where $p_M:S^{*}M\longrightarrow M$ and $p_Y:S^{*}Y\longrightarrow Y,$ are
compatible in the sense that
\[
\gamma (\left( \pi_0) _{!}\left. \sigma _M\right| _{\partial \overline{U}}\right) \alpha
^{-1}=\sigma _Y|_{\partial \overline{Y}}.
\]

Let  $\limfunc{Ell}^k\left( M,\pi \right) $ be the group of stable homotopy
classes of admissible elliptic operators of order $k$ on $M$, modulo elliptic
operators with principal symbols independent of the cotangent variables.

\begin{remark}\label{rem3}
\emph{On manifolds with boundary, one can also consider a similar class of
elliptic operators and boundary value problems. More precisely, let $M$ be a
manifold with boundary, with a projection $\pi$ defined on a closed subset
$\overline{U}\subset M$ as above. We assume that $\overline{U}$ is a
codimension zero submanifold in the interior $M\setminus
\partial M$ and is the Cartesian product $[0,\varepsilon)\times
\overline{U_0}$ in some collar neighbourhood of the boundary for some
codimension zero submanifold $\overline{U_0}$ in $\partial M$. Then on $M$ we
consider operators similar to (\ref{omo1}), where both $D_1$ and $D_2$ are of
order $m$ and are differential operators with respect to the normal variables
in neighbourhoods of the boundaries of the corresponding manifolds (see
(\ref{alef1})). Such operators are considered in the spaces $C^\infty(M,{\cal
E})$. One considers boundary value problems of the form
$$
(D,Bj):C^\infty(M,{\cal E})\longrightarrow C^\infty(M,{\cal F}) \oplus C^\infty(\partial M,{\cal
G}),
$$
where $\mathcal{E},\mathcal{F}\in{\rm Vect}(M,\pi), \mathcal{G}\in{\rm
Vect}(\partial M,\pi|_{\partial M}),$ $j$ is the jet operator of order $m$,
$j:C^\infty(M,{\cal E})\to C^\infty(\partial M,{\cal E}^m|_{\partial M})$, and
the boundary conditions are defined by an admissible operator $B$ on the
boundary.  One can readily extend all results of this section, including the
definition of trivial problems $\cal D_\pm$, the group of stable homotopy
classes of boundary value problems, and order reduction, to this class of
boundary value problems.}
\end{remark}

\begin{remark}
{\em Let $M$ be a manifold with covering $\pi$ on the boundary. In Subsec.~1,
we defined the group ${\rm Ell}^m(M,\pi)$ generated by elliptic non-local
problems for the usual operators. At the same time, the projection $\pi\times
1:\partial M\times [0,1)\to X\times [0,1)$ is defined in a collar neighbourhood
of the boundary, and one can consider the corresponding group ${\rm
Ell}^m(M,\pi\times 1)$ generated by non-local problems for admissible operators
in the sense of Remark~\ref{rem3}. It turns out that these two groups are
isomorphic under the natural mapping
$$
{\rm Ell}^m(M,\pi)\longrightarrow {\rm Ell}^m(M,\pi\times 1).
$$
This essentially follows from the isomorphism
 ${\rm Vect}(M)\simeq{\rm Vect}(M,\pi\times 1).$ }
\end{remark}

\textbf{3. Reduction to a closed manifold.} We return to the problem of
computing the index of non-local operators on a manifold $M$ with a covering
$\pi$ defined on $\partial M$. Consider an embedding  $f:M\rightarrow M^{\prime
}$ in a closed manifold of the same dimension as $M$ (for example, $M^{\prime
}$ can be the double $2M=M\cup _{\partial M}M$). Just as in the classical case
\cite{AtSi1}, $f$ induces the direct image mapping
\[
f_{!}:\limfunc{Ell}^0\left( M,\pi \right) \longrightarrow \limfunc{Ell}^0\left(M^{\prime },\pi
\times 1\right) ,
\]
where $\pi \times 1$ is the extension of  $\pi $ to  $\partial M\times \left[
-1,1\right] \subset M^{\prime }.$ This mapping takes the symbol $\sigma \left(
D\right) =\left(\sigma _M,\sigma _X\right) $ of  an elliptic
operator\footnote{An arbitrary operator $D':C^\infty \left( M,E\right)
\rightarrow C^\infty\left( M,F\right) $ is reduced to this form by adding the
identity operator in the sections of the complementary bundle to $F$.}
\[
D:C^\infty \left( M,E\right) \longrightarrow C^\infty \left( M,\Bbb{C}%
^k\right)
\]
to the symbol on $M^{\prime }$ that coincides on $M$ with the original symbol
and is the identity $id:\Bbb{C}^k\rightarrow \Bbb{C}^k$ on the complement
$M^{\prime }\backslash M$. The extended symbol is defined on the bundle
obtained by clutching $E$ with $\Bbb{C}^k$ using the isomorphism $\left. \sigma
_X\right| _X$ and maps this bundle to the bundle $\Bbb{C}^k$ over the ambient
closed manifold $M^{\prime }.$

\begin{lemma}
The mapping $f_{!}:\limfunc{Ell}^0\left( M,\pi \right) \longrightarrow
\limfunc{Ell}^0\left( M^{\prime },\pi \times 1\right) $ is well defined and is
index preserving.
\end{lemma}

\noindent \emph{Proof}. This is a restatement of the well-known excision
property of the index. The proof is standard, and hence we omit it altogether.
\dokaend

\textbf{4. Embedding in a universal space}. In the index theorems of the
present paper, we assume that the following condition is satisfied.
\begin{assumption}\label{as1}
\emph{The covering $\pi $ is regular and there is a free action of a finite
group $G$ on the submanifold $\overline{U}$ such that $\pi$ is the projection
onto the quotient.}
\end{assumption}

Let  $\left(M,\pi\right)$ and $\left(M^{\prime},\pi^{\prime}\right) $ be two
pairs (both manifolds are assumed to be closed) and let $U$ and $U'$ be the
domains of $\pi$ and $\pi'$, respectively.

\begin{definition}
\emph{We say that $f$ is an } embedding \emph{ of }$\left( M,\pi \right)
$\emph{\ in }$\left( M^{\prime },\pi ^{\prime }\right) $\emph{\ if there is an
embedding }$f:M\rightarrow M^{\prime },$\emph{\ }$f\left( \overline{U}\right)
\subset \overline{U}^{\prime },$\emph{\ that is equivariant on the domain of
$\pi$.}
\end{definition}

Denote by $\pi _N:EG_N\longrightarrow BG_N$ the $N$-universal bundle for $G.$
We assume that $EG_N$ and $BG_N$ are closed manifolds. There is an explicit
construction for such a model (e.g., see \cite{LuMi1}). For example, consider
the embedding
\[
G\subset S_{\left| G\right| }\subset \Bbb{U}\left( \left| G\right| \right)
\]
in the unitary group. (Here $\left| G\right|$ is the order of $G$.) Consider
the bundle $V_{k,\left| G\right| }\rightarrow V_{k,\left| G\right| }/G$, where
$V_{k,n}$ is the Stiefel manifold of  $n$-frames in  $\Bbb{C}^k$. For
sufficiently large $k$, this bundle is $N$-universal.

\begin{proposition}\label{pre1}
For  $\left( M,\pi \right) $ satisfying Assumption \emph{\ref{as1}}, there
exists an embedding in $\left( EG_N,\pi _N\right)$ provided that $N$ is
sufficiently large.
\end{proposition}

\noindent \emph{Proof}. By $N$-universality of $\pi _N$, there exists an
equivariant mapping $\overline{U}\rightarrow EG_N.$ We can assume that this
mapping is a smooth embedding. This can be achieved by a small deformation
provided that the dimension of $EG_N$ is sufficiently large.

This embedding can be extended to a smooth mapping $M\rightarrow EG_N$ owing to
the $N$-connectedness of $EG_N$. Finally, a small deformation outside a
neighbourhood of $U$ makes it a global embedding. \dokaend

\textbf{5. The Euler operator on the disc}. Consider the Neumann problem
\[
\left\{
\begin{array}{c}
\left( d+\delta \right) u=f, \\
\left. \left( \ast u\right) \right| _{\Bbb{S}^{n-1}}=g,
\end{array}
\right. \quad u\in \Lambda ^e\left( \mathbb{D}^n\right) ,f\in \Lambda ^o\left(
\mathbb{D}^n\right),
\]
for the Euler operator in the unit disc $\mathbb{D}^n\subset \Bbb{R}^n$ with
the Euclidean metric. Here $g\in \Lambda ^{e+n}\left( \Bbb{S}^{n-1}\right) .$
This boundary value problem is elliptic, and Hodge theory shows that the
cokernel is trivial and the one-dimensional kernel consists of constant
functions.

The same is true for the homogeneous boundary value problem, which we rewrite in the operator
form
\[
D_{dR}=d+\delta :\Lambda _0^e\left( \mathbb{D}^n\right) \longrightarrow \Lambda ^o\left(
\mathbb{D}^n\right).
\]
(Here $\Lambda _0^e\left( \mathbb{D}^n\right) $ is the space of forms
satisfying the homogeneous boundary condition.) This operator is $O\left(
n\right) $-equivariant with respect to the natural action of the orthogonal
group on $\mathbb{D}^n.$\vspace{1mm}

\textbf{6. Embeddings and the index of elliptic operators}. Let $f:\left( M,\pi
\right)\rightarrow \left( M^{\prime },\pi ^{\prime }\right)$ be an embedding of
positive codimension. We choose a Riemannian metric on  $M^{\prime }$ that is
$G$-invariant over $\overline{U}^{\prime }\subset M^{\prime }$. Denote the
normal bundle to $M$ by $NM.$ Then a closed tubular neighbourhood   $W$ of  $M$
in $M^{\prime } $ is diffeomorphic to the unit ball subbundle $DM\subset NM.$
Additionally, we can assume this diffeomorphism to be $G$-equivariant over
$\overline{U}\subset \overline{U}^{\prime }$.

Consider an admissible elliptic operator
\[
D:C^\infty \left( M,{\cal E}\right) \longrightarrow C ^\infty \left( M,{\cal F}\right) .
\]

We define a boundary value problem on $DM$ as the exterior tensor product of
$D$ by a family of boundary value problems for the Euler operator in the
fibers. The definition of this product is the same as in \cite{AtSi1}.

More precisely, the exterior tensor product gives the operator
\[
\mathcal{D}=\left(
\begin{array}{cc}
\widetilde{D}\otimes 1_{\Lambda ^e} & -1_F\otimes \widetilde{D}_{dR}^{*} \\
1_E\otimes \widetilde{D}_{dR} & \widetilde{D}^{*}\otimes 1_{\Lambda ^o}
\end{array}
\right),
\]
where the pullback of $D$ to the bundle $DM$ with coefficients in even forms on
the fibers is denoted by
$$
\widetilde{D}\otimes 1_{\Lambda ^e}: C^\infty \left( DM,p^{*}{\cal E}\otimes \Lambda _0^e\left(
DM\right) \right) \longrightarrow C^\infty \left( DM,p^{*}{\cal F}\otimes \Lambda _0^e\left(
DM\right) \right).
$$
Here $\widetilde{D}^{*}\otimes 1_{\Lambda ^o}$ is the pullback of the adjoint
operator with coefficients in odd forms. The family of Neumann problems for the
Euler operator $\widetilde{D}_{dR}^M$ with coefficients in the triple ${\cal
E}=\left(E,E_0,\alpha \right)$ is denoted by
$$
1_E\otimes \widetilde{D}_{dR}:C_{\alpha \otimes 1}^\infty \left( DM,p^{*}E\otimes \Lambda
_0^e\left( DM\right) \right) \longrightarrow C_{\alpha \otimes 1}^\infty \left( DM,p^{*}E\otimes
\Lambda ^o\left( DM\right) \right).
$$
The off-diagonal entries of $\mathcal{D}$ commute with entries on the diagonal
by construction. As in ordinary Atiyah--Singer theory, this leads to the
following result.

\begin{lemma} One has
$\limfunc{ind}D=\limfunc{ind}\mathcal{D}.$
\end{lemma}

\noindent The \emph{proof} is similar to \cite{AtSi1}. \dokaend

Thus an elliptic operator on the submanifold $M\subset M^{\prime }$ induces an
elliptic boundary value problem with the same index on the tubular
neighbourhood $DM\simeq W\subset M^{\prime }$. Further, we can apply the order
reduction procedure to this problem  (see Remark~\ref{remk}) and extend the
resulting zero-order operator from  $W$ to the entire manifold $M^{\prime} $ as
in Subsec.~2.

Summarizing, we see that the embedding
 $f$ of $\left(M,\pi\right)$ in
$\left(M^{\prime},\pi^{\prime}\right) $ induces the direct image mapping
\[
f_{!}:\limfunc{Ell}^1\left( M,\pi \right) \longrightarrow \limfunc{Ell}^0\left( M^{\prime },\pi
^{\prime }\right),
\]
which preserves the index.

\begin{remark}
\emph{A straightforward computation shows that the linear homotopy of order
reduction for boundary value problems (defined in Remark~\ref{remk}) which
extends the symbol }$\sigma \left( d+\delta \right) $\emph{\ from
}$T^{*}\mathbb{D}^n$\emph{\ to }$T^{*}\Bbb{R}^n$\emph{\ as an invertible
element outside a compact set defines an element of the equivariant
}$K$-\emph{group equal to the element}
\[
j_{!}\left( 1\right) \in K_{O\left( n\right) }\left( T^{*}\Bbb{R}^n\right) ,\quad
j:pt\longrightarrow \Bbb{R}^n,
\]
\emph{\ which is used in the standard proof of the Atiyah--Singer theorem.}
\end{remark}

\textbf{7. The Index theorem}. Let $f$ be an embedding of $(M,\pi)$ in the
universal space defined in Proposition~\ref{pre1}. For the universal space
$EG_N$, the projection $\pi _N$ is defined globally. Therefore, the direct
image of a non-local operator can be treated as a usual elliptic operator on
the base  $BG_N$; i.e., we have a natural mapping
\[
\left( \pi _N\right) _{!}:\limfunc{Ell}\left( EG_N,\pi _N\right) \longrightarrow
\limfunc{Ell}\left( BG_N\right) \simeq K\left( T^{*}BG_N\right) .
\]

\begin{theorem}\label{thind1}
For a pair $(M,\pi)$ satisfying Assumption~\emph{\ref{as1}}, the diagram
\[
\begin{array}{ccc}
\limfunc{Ell}^1\left( M,\pi \right)  & \stackrel{f_{!}}{\longrightarrow } &
\limfunc{Ell}^0\left( EG_N,\pi _N\right)  \\
\limfunc{ind}\downarrow \quad  &  & \downarrow \left( \pi _N\right) _{!} \\
\quad\Bbb{Z} & \stackrel{\limfunc{ind}_t}{\longleftarrow } & K\left( T^{*}BG_N\right) ,
\end{array}
\]
commutes. Here $\limfunc{ind}_t$  is the usual topological index on a closed manifold.
\end{theorem}

\noindent \emph{Proof}. Indeed, we have
\[
\limfunc{ind}D=\limfunc{ind}f_{!}\left[ D\right] =\limfunc{ind}\left( \pi _N\right)
_{!}f_{!}\left[ D\right] =\limfunc{ind}_t\left( \left( \pi _N\right) _{!}f_{!}\left[ D\right]
\right).
\]
The first equality here follows from the invariance of the index for
embeddings, the second from the fact that $\left( \pi_N\right) _{!}$ does not
change the operator, and the last equality is just the Atiyah--Singer formula
on $BG_N.$ \dokaend

\textbf{8. Example}. \emph{Manifolds with reflecting boundary} \cite{Hsi2}.\
Let $M$\ be a $4k$-dimensional compact oriented Riemannian manifold with
boundary $\partial M$. Suppose that $\partial M$ is equipped with an
orientation-reversing smooth involution $G$ without fixed points. The
involution defines a free action of the group $\Bbb{Z}_2$\ and the
corresponding double covering $\pi :\partial M\longrightarrow \partial
M/\Bbb{Z}_2.$ Consider the Hirzebruch operator \cite{Pal1}
\[
d+d^{*}:\Lambda ^{+}\left( M\right) \longrightarrow \Lambda ^{-}\left( M\right).
\]
In a neighbourhood of the boundary, let us take a metric lifted from $\left[
0,1\right] \times
\partial M/\Bbb{Z}_2.$ Then the Hirzebruch operator can be decomposed near the boundary  as (see \cite{APS1})
\[
\frac \partial {\partial t}+A
\]
(up to a bundle isomorphism), where $A$ is an elliptic self-adjoint operator on
the boundary and is given by the formula
\[
A:\Lambda ^{*}\left( \partial M\right) \longrightarrow \Lambda ^{*}\left(
\partial M\right) ,\qquad
A\omega =\left( -1\right) ^{k+p}\left( d*-\varepsilon *d\right) \omega;
\]
here for an even degree form $\omega \in \Lambda ^{2p}\left(
\partial M\right) $\ we set $\varepsilon =1,$\ and $\varepsilon
=-1$ otherwise. Since $G$ reverses the orientation, it follows that $A$ and
$G^*$ anticommute:
\[
G^{*}A=-AG^{*}.
\]
It is known that the Hirzebruch operator has no well-posed classical boundary
conditions. However, it admits the non-local boundary value problem
\begin{equation}
\left\{
\begin{array}{ll}
\left( d+d^{*}\right) \omega =f,\vspace{2mm} &  \\
\frac{\left( 1+G^{*}\right) }2\left. \omega \right| _{\partial M}=g, & \quad g\in \Lambda
^{*}\left( \partial M\right) ^{\Bbb{Z}_2}\simeq \Lambda ^{*}\left( \partial M/\Bbb{Z}_2\right)
\end{array}
\right.  \label{nel1}
\end{equation}
on the manifold with reflecting boundary. Here $\Lambda ^{*}\left( \partial
M\right) ^{\Bbb{Z}_2}$ is the subspace of $G$-invariant forms on the boundary.

\begin{proposition}
The non-local boundary value problem \emph{(\ref{nel1})} is elliptic.
\end{proposition}

\emph{Proof}. Consider an arbitrary point $x\in \partial M/\Bbb{Z}_2$. An
explicit computation shows that near this point the equivalent classical
boundary value problem is
\[
\left\{
\begin{array}{c}
\left( \frac \partial {\partial t}+A\right) \omega _1=f_1,\quad \left( \frac
\partial {\partial t}-A\right) \omega _2=f_2,\vspace{2mm} \\
\left. \omega _1\right| _{\partial M/\Bbb{Z}_2}+\left. \omega _2\right| _{\partial
M/\Bbb{Z}_2}=g.
\end{array}
\right.
\]
It is elliptic (satisfies the Shapiro--Lopatinskii condition), since the symbol
of the operator of boundary conditions defines an isomorphism
\[
\func{Im}\sigma \left( \Pi _{+}\right) \left( x,\xi \right) \oplus \func{Im}%
\sigma \left( \Pi _{-}\right) \left( x,\xi \right) \simeq \Lambda ^{*}\left(
\partial M\right) _x
\]
at an arbitrary point $\left( x,\xi \right) \in S^{*}(\partial M/\Bbb{Z}_2)$,
where
\[
\Pi _{+}=\frac{A+\left| A\right| }{2\left| A\right| }
\]
is the non-negative spectral projection of $A$ and $\Pi _{-}=1-\Pi _{+}$ is the
negative projection. The ellipticity of the classical boundary value problem
proves the desired statement. \dokaend

\begin{proposition}
One has
\[
\limfunc{ind}\left( d+d^{*},{\left( 1+G^{*}\right) }\right) =\limfunc{%
sign}M,
\]
where $\limfunc{sign}M$ is the signature of $M$.
\end{proposition}

\emph{Proof}. The symbol of (\ref{nel1}) coincides with that of the composition
of the spectral Atiyah--Patodi--Singer boundary value problem
\[
\left\{
\begin{array}{ll}
\left( d+d^{*}\right) \omega =f &  \\
\Pi _{+}\left. \omega \right| _{\partial M}=\omega ^{\prime }, & \quad \omega ^{\prime }\in
\func{Im}\Pi _{+}\subset \Lambda ^{*}\left( \partial M\right) ,
\end{array}
\right.
\]
and the Fredholm operator
\begin{equation}
{\left( 1+G^{*}\right) }:\func{Im}\Pi _{+}\longrightarrow \Lambda ^{*}\left( \partial M\right)
^{\Bbb{Z}_2}.  \label{subs4}
\end{equation}
Let us compute both indices.

1) For the index of the spectral boundary value problem, one has \cite{APS1}
\[
\limfunc{ind}\left( d+d^{*},\Pi _{+}\right) =\limfunc{sign}M-\frac{\dim \ker A}2.
\]
In addition, by the Hodge--de Rham theory we obtain $ {\dim \ker A}={\dim
H^{*}\left(
\partial M\right) } $.

2) On the other hand, one can readily verify that the operator in
Eq.~(\ref{subs4}) is surjective and its kernel coincides with the space of
$G$-antiinvariant harmonic forms. The Hodge operator $*$ interchanges the
antiinvariant and invariant subspaces. Thus we obtain
\[
\dim \ker \left. \left( 1+G^{*}\right) \right| _{\func{Im}\Pi _{+}\left( A\right) }=\frac{\dim
\ker A}2.
\]
Adding the index of the spectral problem to the index of $\left( 1+G^{*}\right)
,$ we obtain
\[
\limfunc{ind}\left( d+d^{*},{\left( 1+G^{*}\right) }\right) =\limfunc{%
sign}M-\frac{\dim \ker A}2+\frac{\dim \ker A}2=\limfunc{sign}M.
\]

The proof of the theorem is complete. \dokaend

\section{The homotopy classification of non-local operators}\label{fama1}

Let us cut  $M$ into two parts
\[
M^{\prime }=M\backslash \left\{ \partial M\times \left[ 0,1\right) \right\}
\simeq M\;\qquad \text{and\qquad }\partial M\times \left[ 0,1\right] .
\]
Then the symbol $\sigma \left(D\right) $ of an admissible elliptic operator $D$
of order zero is naturally represented as a pair of usual symbols
\begin{equation}
\left. \sigma \left( D\right) \right| _{M^{\prime }}\text{ and }\left( \pi
\times 1\right) _{!}\left. \sigma \left( D\right) \right| _{\partial M\times
\left[ 0,1\right] }.  \label{smb1}
\end{equation}
Both symbols define  difference elements
\[
\left[ \sigma _{M}\right] \in K\left( T^{*}M^{\prime }\right) ,\qquad \left[ \sigma_X\right] \in
K\left( T^{*}\left( X\times \left( 0,1\right] \right) \right).
\]
Here and in what follows, we use $K$-groups with compact supports. In the
latter case, the elliptic symbol
 $\sigma_X $ of order zero is invertible over $X\times \left\{
0\right\} $ (this follows from ellipticity and the decomposition in
Eq.~(\ref{alef1})) and hence defines element in the above-mentioned $K$-group
with compact supports.

However, it is impossible to define an element of a single topological
$K$-group; indeed, the manifolds  $T^{*}M^{\prime }$ and $T^{*}\left( X\times
\left( 0,1\right] \right) $ can not be glued together, for their boundaries are
not diffeomorphic. Nonetheless, we can glue the algebras of functions on these
spaces instead of the original manifolds.

\textbf{1. The $C^*$-algebra of a manifold with a covering on the boundary.} To
each space, we assign an algebra of continuous functions vanishing at infinity:
\[
C_0\left( T^{*}M^{\prime }\right),\quad C_0\left( T^{*}\left( X\times \left( 0,1\right] \right)
,\limfunc{End}p^{*}\pi _{!}1\right).
\]
More precisely, on the space $T^*(X\times (0,1])$ we consider functions ranging
in the set of endomorphisms of the bundle $\pi _{!}1\in \limfunc{Vect}\left(
X\right) $, where $p:T^{*}\left( X\times \left( 0,1\right] \right) \rightarrow
X$ is the natural projection. In the direct sum of these algebras, consider the
subalgebra determined by the compatibility condition
\begin{equation}
\label{inty5} \mathcal{A}_{T^{*}M,\pi }=\left\{ (u,v)\;\left|
\begin{array}{c}
u\in C_0\left( T^{*}M^{\prime }\right), v\in C_0\left( T^{*}\left( X\times \left( 0,1\right]
\right) , \limfunc{End}p^{*}\pi _{!}1\right)
\vspace{1mm}\\
\;\beta \left. u\right| _{\partial M^{\prime }}\beta ^{-1}=\left. v\right| _{t=1}
\end{array}
\right.\right\}.
\end{equation}
Here $t$ is the coordinate on $\left( 0,1\right] $.

{For the trivial covering }$\partial M\rightarrow \partial M=X${, this algebra
is just the commutative algebra of continuous functions on }$T^{*}\left(
M\backslash
\partial M\right)$ {\ vanishing at infinity. Let us also mention that this algebra can be also
viewed as the groupoid $C^*$-algebra \cite{BrNi1}  of the equivalence relation
defined by $\pi$ ($x\sim y$ if  either $x=y$ or  $x,y\in \partial M$ and
$\pi(x)=\pi(y)$). For a trivial covering, this algebra was used in
\cite{Ros2}.}

\textbf{2. The difference construction.} Let us define the \emph{difference
construction} for non-local operators. This will be a mapping
\begin{equation}
\chi :\limfunc{Ell}^0\left( M,\pi \right)
\longrightarrow K_0\left( %
\mathcal{A}_{T^{*}M,\pi }\right)  \label{stk}
\end{equation}
into the $K_0$  group of the $C^{*}$-algebra $\mathcal{A}_{T^{*}M,\pi }.$ To
this end, we take an elliptic operator
\[
D:C^\infty \left( M,E\right) \longrightarrow C^\infty \left(
M,{\Bbb{C}^k}\right),
\]
fix some embeddings of  $E$ and ${\Bbb{C}^k}$ in trivial bundles of
sufficiently large dimension, and denote by  $P_E$ and $P_{\Bbb{C}^k}$ the
projections that define the corresponding subbundles:
\[
E\simeq \func{Im}P_E\subset \Bbb{C}^N\oplus 0,\qquad {\Bbb{C}^k}\simeq \func{Im}%
P_{\Bbb{C}^k}\subset 0\oplus \Bbb{C}^L.
\]
We denote the direct images of these projections near the boundary by  $P_{\pi
_{!}E}$ and $P_{\pi _{!}{\Bbb{C}^k}}$.  The  \emph{difference element of $D$}
is, by definition, the difference
\[
\chi \left[ D\right] =\left[ P_1\oplus P_2\right] -\left[ P_{\Bbb{C}^k}\oplus P_{\pi%
_{!}{\Bbb{C}^k}}\right]\in K_0\left( %
\mathcal{A}_{T^{*}M,\pi }\right) ,
\]
where the projection  $P_1$ over $M^{\prime }$ is given by
\begin{equation}
\left\{
\begin{array}{cc}
P_E\cos ^2\left| \xi \right| + P_{\Bbb{C}^k}\sin ^2\left| \xi \right| + \left(\sigma^{-1}_M
\left( x,\xi \right)P_{\Bbb{C}^k} +\sigma_M \left( x,\xi \right) P_E \right) \sin \left| \xi
\right|  \cos \left| \xi \right|
 , & \left| \xi
\right| \leq \pi /2,\vspace{2mm} \\
P_{\Bbb{C}^k}, & \left| \xi \right| >\pi /2.
\end{array}
\right.  \label{str}
\end{equation}
(We assume that the principal symbol is zero-order homogeneous in  $\xi$.) The
projection $P_2$ over $X\times \left[ 0,1\right] $ is defined by the formula
\[
P_2=
\left\{
\begin{array}{cc}
P_{\pi _{!}E}\cos ^2\left| \xi \right| +P_{\pi _{!}{\Bbb{C}^k}}\sin ^2\left| \xi%
\right| +1/2 \left(\widetilde{\sigma }^{-1}\left( x^{\prime },\xi \right) P_{\pi
_{!}{\Bbb{C}^k}}+ \widetilde{\sigma }\left( x^{\prime },\xi \right) P_{\pi _{!}E} \right)\sin
2\left| \xi \right|
, & \vspace{2mm} \\
P_{\pi _{!}E}\cos ^2\varphi +P_{\pi _{!}{\Bbb{C}^k}}\sin ^2\varphi + 1/2\left(\widetilde{\sigma
}^{-1}\left( x^{\prime },0\right) P_{\pi _{!}{\Bbb{C}^k}}+ \widetilde{\sigma }\left( x^{\prime
},0\right) P_{\pi _{!}E} \right) \sin 2\varphi, &
\vspace{2mm} \\
P_{\pi _{!}{\Bbb{C}^k}}, &
\end{array}
\right.
\]
where the first case is used for $x^{\prime }\in X\times \left[ 1/2,1\right]
,\left| \xi \right| \leq \pi /2 $, the second for $x^{\prime }\in X\times
\left[ 0,1/2\right] ,\left| \xi \right| <\pi t$, and the third otherwise. Here
we write
\[
\varphi =\left| \xi \right| +\pi /2\left( 1-2t\right) ,\quad \widetilde{%
\sigma }\left( x^{\prime },\xi \right) =\sigma_X  \left( x^{\prime },\xi
\right)
\]
for brevity. Geometrically, these projections define a subbundle that coincides
with $E\subset \Bbb{C}^{N+L}$ over the zero section (for $\xi$=0); coincides
with the orthogonal bundle ${\Bbb{C}^k}\subset \Bbb{C}^{N+L}$ for $\left| \xi
\right| \geq \pi /2$;  and is  obtained by the rotation of the first bundle
towards the second bundle with the use of $\sigma \left( D\right) $ at the
intermediate points. (The symbol is treated as an isomorphism of the two
bundles.) By construction, $P_1$ and $P_{\Bbb{C}^k}$ coincide outside a compact
set in $T^{*}M^{\prime },$ and $P_2$ and $P_{\pi _{!}{\Bbb{C}^k}}$ coincide
outside a compact set in $T^{*}\left( X\times \left( 0,1\right] \right).$
Therefore, the difference $\left[ P_1\oplus P_2\right] -\left[
P_{\Bbb{C}^k}\oplus P_{\pi _{!}{\Bbb{C}^k}}\right] $ is indeed in $K_0\left(
\mathcal{A}_{T^{*}M,\pi }\right) .$

\begin{remark}
\emph{This element of the }$K$\emph{-group can be equivalently defined by
different expressions
(cf. \cite{Hig2}).}%
\end{remark}

\begin{theorem}
\label{sem1}The difference construction $\emph{(\ref{stk})}$ is a
well-defined group isomorphism.
\end{theorem}

\emph{Proof.} The mapping $\chi $ preserves the equivalence relations in
$\limfunc{Ell}^0\left(M,\pi \right) $ and $K_0\left( \mathcal{A}_{T^{*}M,\pi
}\right) .$ Indeed, under an operator homotopy the symbols vary continuously.
Therefore, the corresponding projections $P_{1,2}$ are joined by a continuous
homotopy. Furthermore, $\chi \left[ D\right] $ is independent of the choice of
an embedding in a trivial bundle, since all such embeddings are homotopic, and
for a trivial $D$ (i.e., one induced by a vector bundle isomorphism) $\chi[D]$
is equal to zero. This shows that $\chi$ is well defined. The proof that this
mapping is one-to-one presents no essential difficulties and is left to the
reader.

\dokaend\vspace{1mm}

\textbf{3. Index theorem for families}. Later on in Section \ref{defe3}, we use
a families index formula. Let us briefly state the corresponding results.

Let $P$ be a compact space. Denote by $\limfunc{Ell}_P\left( M,\pi \right) $ the group of stable
homotopy classes of elliptic families on $M$ parametrized by $P.$

\begin{theorem}[\rm the index of families of non-local operators]\label{indfam} Let $(M,\pi)$ satisfy
Assumption~{\rm \ref{as1}}. Then for an embedding $f:M\to EG_N$ the direct
image mapping $f_!$ for families is well defined and the following diagram
commutes:
$$
\xymatrix{ \limfunc{Ell}_P\left( M,\pi \right)\ar[d]_{\rm ind} \ar[rr]^{f_!}
& & \limfunc{Ell}_P\left( EG_N,\pi_N \right) \ar[d]^{(\pi_N)_!}\\
K^0(P) & \ar[l]_{{\rm ind}_t \qquad} K^0\left(P\times T^{*} BG_N\right) \ar@{=}[r]^{\;} &
\limfunc{Ell}_P\left( BG_N \right) . }
$$
\end{theorem}
The proof is similar to that of Theorem~\ref{thind1} in the previous section
(cf. \cite{AtSi4}) and therefore is omitted.

Let us finally note that the difference construction can also be defined in
this case as a mapping
\[
\chi _P:\limfunc{Ell}_P\left( M,\pi \right) \longrightarrow%
K_0\left( C\left( P,\mathcal{A}_{T^{*}M,\pi }\right) \right) ,
\]
where $C\left( P,\mathcal{A}_{T^{*}M,\pi }\right) $ is the algebra of
continuous functions on $P$ ranging in the $C^{*}$-algebra
$\mathcal{A}_{T^{*}M,\pi }.$

\section{A homotopy invariant for manifolds with covering on the boundary
\label{pr1}}

\textbf{1. The class of operators}. On a manifold $M$ with covering $\pi $ on
the boundary, we consider elliptic differential operators
\[
D:C^\infty \left( M,E\right) \longrightarrow C^\infty \left( M,F\right)
\]
that are lifted from the base of the covering in a neighbourhood of the
boundary . Technically, we suppose that the following condition is satisfied.

\begin{assumption}\label{as2}{\em
The restrictions of the bundles $E$ and $F$ to the boundary are lifted from the
base of the covering; moreover, we fix some isomorphisms
\[
\left. E\right| _{\partial M}\simeq \pi ^{*}E_0,\quad \left. F\right| _{\partial M}\simeq \pi
^{*}F_0,\quad \quad E_0,F_0\in \limfunc{Vect}\left( X\right) ,
\]
and for some operator  $ D_0:C^\infty \left( X\times \left[ 0,1\right)
,E_0\right) \longrightarrow C^\infty \left( X\times \left[ 0,1\right)
,F_0\right)$ on the cylinder with base $X$ the direct image of $D$ in a collar
neighbourhood of the boundary satisfies the commutative diagram
\begin{equation}
\xymatrix{ C^\infty \left( X\times \left[ 0,1\right) ,\pi _{!}E\right) \ar[d]_\simeq
\ar[r]^{\left( \pi \times 1\right) _{!}D}& C^\infty \left(X\times \left[ 0,1\right) ,\pi
_{!}F\right) \ar[d]^\simeq
\\
C^\infty \left( X\times \left[ 0,1\right) ,E_0\otimes \pi _{!}1\right) \ar[r]^{D_0\otimes 1} &
C^\infty \left( X\times \left[0,1\right) ,F_0\otimes \pi _{!}1\right) . } \label{pdown}
\end{equation}
Here $D_0\otimes 1$ stands for the operator $D_0$ with coefficients in the flat
bundle $\pi _{!}1$ (e.g., see \cite{APS3}).}
\end{assumption}

We also suppose that $D$ is first-order operator and the following assumption
is satisfied.

\begin{assumption}\label{as3}{\em
In the neighbourhood $X\times\left[ 0,\varepsilon \right) $ of the boundary,
the operator has the form
\[
\left. D_0\right| _{X\times \left[ 0,\varepsilon \right) }=\Gamma \left( \frac \partial {\partial
t}+A_0\right)
\]
for a bundle isomorphism $\Gamma $, where $A_0$ is an elliptic self-adjoint first-order operator
on $X.$ This operator is  called the \emph{tangential operator} of $D_0.$}
\end{assumption}

If $D$ satisfies Assumptions \ref{as2} and \ref{as3},  then near the boundary it has the form
\[
D=\frac \partial {\partial t}+\pi ^{!}\left( A_0\otimes 1\right)
\]
up to a vector bundle isomorphism.  For brevity, the self-adjoint operator $\pi
^{!}\left( A_0\otimes 1\right) $ will be denoted by $A.$\vspace{2mm}

\textbf{2. The homotopy invariant}. For an operator $D$ satisfying
Assumptions~\ref{as2} and \ref{as3}, consider the \emph{spectral
Atiyah--Patodi--Singer boundary value problem} \cite{APS1}
\[
\left\{
\begin{array}{llll}
Du & = & f, &  \\
\Pi _{+}\left. u\right| _{\partial M} & = & g, & \qquad g\in \func{Im}\Pi
_{+}\subset H^{s-1/2}\left( \partial M,E\right) ,
\end{array}
\right.
\]
where $ \Pi _{+}=(A+\left| A\right| )/{2\left| A\right| }$ is the non-negative
spectral projection of the self-adjoint operator  $A$. (If $A$ is not
invertible, then in this formula  one should replace $A$ by $A+\varepsilon$ for
some $\varepsilon$ less then the absolute value of the greatest negative
eigenvalue of $A$.) The spectral problem is always Fredholm. However, its index
$\limfunc{ind}\left( D,\Pi _{+}\right) $ is not invariant under homotopies of
$D$ and is not determined by its principal symbol. Here by definition a
continuous homotopy of $D$ is a continuous homotopy in the interior of $M$ that
can be covered in a neighbourhood of the boundary by a homotopy of the diagram
\eqref{pdown}  and a continuous homotopy of tangential operators.

\begin{proposition}
The sum
\begin{equation}
\widetilde{\limfunc{ind}}D\stackrel{def}{=}\func{mod}n\text{-}\left(
\limfunc{ind}\left( D,\Pi _{+}\right) +\eta \left( A\right) -n\eta
\left( A_0\right) \right) \in \Bbb{R}/n\Bbb{Z},  \label{moda1}
\end{equation}
is a homotopy invariant of $D.$ Here $n$ is the number of sheets of the
covering and $\eta \left( A\right)$ and $\eta \left( A_0\right) $ are the
spectral Atiyah--Patodi--Singer $\eta $-invariants of the tangential operators
$A$ and $A_0$.
\end{proposition}

\emph{Proof}. Consider the non-reduced invariant
\begin{equation}
\limfunc{ind}\left( D,\Pi _{+}\right) +\eta \left( A\right) -n\eta \left(
A_0\right) .  \label{bet3}
\end{equation}
The results of \cite{APS3} imply that for a smooth operator family $D_t$ this
expression is a piecewise smooth function of the parameter $t$ (the
corresponding families of tangential operators are denoted by $A_t$ and
$A_{0,t}$).

1) We claim that (\ref{bet3}) is a piecewise constant function. Indeed, the derivative of the
$\eta $-invariant
\[
\frac d{dt}\eta \left( A_t\right)
\]
with respect to $t$ is local, i.e., is equal to an integral over the manifold
of an expression determined by the complete symbol of the tangential family
$A_t.$ However, the complete symbols of $A_t$ and $A_{0,t}$ coincide locally by
Assumption \ref{as2}. Thus we have
\[
\frac d{dt}\eta \left( A_t\right) =n\frac d{dt}\eta \left( A_{0,t}\right) .
\]
Therefore, (\ref{bet3}) is a piecewise constant function.

2) Let us show that the jumps of this function are multiples of the number of
sheets of the covering. Indeed, for a homotopy $D_t$ the index, as well as the
$\eta $-invariants, changes by the spectral flow of the corresponding families
of tangential operators. Hence
\[
\left. \left[ \limfunc{ind}\left( D_t,\Pi _{+,t}\right) +\eta \left(
A_t\right) -n\eta \left( A_{0.t}\right) \right] \right| _{t=0,1}=\left(
-1+1\right) \limfunc{sf}\left( A_t\right) _{t\in \left[ 0,1\right] }-n%
\limfunc{sf}\left( A_{0,t}\right) _{t\in \left[ 0,1\right] }\in n\Bbb{Z},
\]
as desired. \dokaend

\begin{remark}
\emph{For a trivial covering, our invariant is none other than the
}$\func{mod}n$-index
\[
\func{mod}n\text{-}\limfunc{ind}(D,\Pi_+)\in \Bbb{Z}_n\subset \Bbb{R}/n\Bbb{Z}
\]\emph{\
of Freed-Melrose} \emph{%
\cite{FrMe1}}. \emph{On the other hand, the fractional part of the invariant}
\emph{\ (\ref{moda1}) is the so-called} relative Atiyah--Patodi--Singer $\eta
$-invariant
\emph{\cite{APS2,APS3}} %
\[
\left\{ \eta \left( A_0\otimes 1_{\pi _{!}1}\right) -n\eta \left( A_0\right)
\right\} \in \Bbb{R}/\Bbb{Z}
\]
\emph{of }$A_0$\emph{\ with coefficients in the flat bundle }$\pi _{!}1\in \limfunc{%
Vect}\left( X\right) .$
\end{remark}

The invariant $\widetilde{\limfunc{ind}}$ has an interesting interpretation as an obstruction.
Namely, suppose that $M$ is the total space of a covering $\widetilde{\pi }$ with base $Y$ that
induces the covering $\pi $ over the boundary
\[
\begin{array}{ccc}
\partial M & \subset & M \\
\pi \downarrow \;\; &  & \;\;\downarrow \widetilde{\pi } \\
X & \subset & Y.
\end{array}
\]

\begin{proposition}
If a differential operator $D:C^\infty \left( M,E\right) \rightarrow C^\infty
\left( M,F\right) $ is the pullback of an elliptic operator $D_0$ on $Y,$ then
\[
\widetilde{\limfunc{ind}}D=0.
\]
\end{proposition}

\emph{Proof}. According to the Atiyah--Patodi--Singer formula (see \cite{APS1}), the sum
\[
\limfunc{ind}\left( D,\Pi _{+}\right) +\eta \left( A\right)
\]
is equal to the integral over the manifold of a local expression defined by the
complete symbol of $D$. Since $D$ and $D_0$ coincide locally, one has
\[
\limfunc{ind}\left( D,\Pi _{+}\right) +\eta \left( A\right) =n\left(
\limfunc{ind}\left( D_0,\Pi _{+,0}\right) +\eta \left( A_0\right) \right) .
\]
We obtain the desired formula by transposing the term $n\eta \left( A_0\right)
$ to the left-hand side. \dokaend

\section{The index defect formula \label{defe3}}

The aim of this section is to find a topological formula for the invariant
$\widetilde{\rm ind}$.\vspace{1mm}

\textbf{1. The difference construction}. The pair $\left( M,\pi \right) $
defines the singular space \fgr{311_2.eps}{fii2}{Singular space}
\[
\overline{M}^\pi =M\left/ \left\{ x\sim x^{\prime },\text{
if }x,x^{\prime }\in \partial M\text{ and }\pi \left(
x\right) =\pi \left( x^{\prime }\right) \right\} \right. ,
\]
obtained by identification of points in the fibers of $\pi$ (see Figure
\ref{fii2} in the case of a trivial covering). Likewise, the boundary of the
non-compact manifold $T^{*}M$ is a covering over the product $T^{*}X\times
\Bbb{R}$, and the corresponding singular space will be denoted by
$\overline{T^{*}M}^\pi$.

Consider an elliptic operator $D$ satisfying Assumption \ref{as2}. The diagram
\eqref{pdown} implies that the principal symbol defines a $K$-theory element
\[
\left[ \sigma \left( D\right) \right] \in K\left( \overline{T^{*}M}^\pi
\right) .
\]
Thus we have a homomorphism
\[
\chi :  \limfunc{Ell}\left( \overline{M}^\pi \right)\longrightarrow K\left( \overline{T^{*}M}^\pi
\right) ,\quad
 \chi \left[ D\right]  =  \left[ \sigma \left( D\right) \right] .
\]
Here $\limfunc{Ell}\left( \overline{M}^\pi \right) $ is the Grothendieck group
of homotopy classes of elliptic operators $D$ on $M$ that satisfy Assumptions
\ref{as2} and \ref{as3}.

The topological formula for the invariant $\widetilde{\limfunc{ind%
}}$ uses the Poincar\'{e} pairing on the manifold
$\overline{T^{*}M}^\pi $ with singularities. Let us define this
pairing. \vspace{2mm}

\textbf{2. A pairing in }$K$\textbf{-theory of a singular manifold}. By analogy
with the algebra $\mathcal{A}_{T^{*}M,\pi }$ of the cotangent bundle, one can
define an algebra for $M$ itself:
\[
\mathcal{A}_{M,\pi }=\left\{ (u,v)\; \left|\;
\begin{array}{c}
u\in C_0\left( M^{\prime }\right), v\in C_0\left(
X\times \left( 0,1\right] ,\limfunc{End}\pi _{!}1\right)\vspace{1mm}\\
\beta \left( \left. u\right| _{\partial M^{\prime }}\right)\beta ^{-1}=\left. v\right| _{t=1}
\end{array}
\right. \right\} .
\]

\begin{lemma}
\label{lemc2}The group $K_0\left( \mathcal{A}_{M,\pi}\right) $ is isomorphic
to the group of stable homotopy classes of triples
\[
\left( E,F,\sigma \right) ,\qquad E,F\in \limfunc{Vect}\left( M\right) ,\;\sigma :\pi _{!}\left.
E\right| _{\partial M}\longrightarrow \pi _{!}\left. F\right| _{\partial M}.
\]
Here $\sigma $ is a bundle isomorphism, and trivial triples are those with
$\sigma $ induced by an isomorphism over $M$.
\end{lemma}

\emph{Proof.} Note that this lemma is similar to  Theorem~\ref{sem1}, which can
be also considered as giving a realization of the group $K_0\left(
\mathcal{A}_{T^{*}M,\pi }\right)$ in topological terms. Along the same lines, a
triple $\left( E,\mathbb{C}^k,\sigma \right) $ defines the element
\[
\left[ P_E\oplus P_2\right] -\left[ P_{\mathbb{C}^k}\oplus P_{\pi _{!}\mathbb{C}^k}\right] \in
K_0\left( \mathcal{A}_{M,\pi }\right) ,
\]
where the projection $P_2$ over $X\times \left[ 0,1\right] $ is
\[
P_2=P_{\pi _{!}E}\cos ^2\varphi +P_{\pi _{!}{\mathbb{C}^k}}\sin ^2\varphi +P_{\pi%
_{!}{\mathbb{C}^k}}\sigma \left( x\right) P_{\pi _{!}E}\sin 2\varphi ,\quad \varphi =\frac \pi
2\left( 1-t\right) .
\]
Here $P_E,P_{\mathbb{C}^k}$ are projections on subbundles isomorphic to $E$ and
$F$. We also suppose that the subbundles are orthogonal to each other.

The proof of the fact that this mapping induces an isomorphism with the group
$K_0\left( \mathcal{A}_{M,\pi }\right) $ is similar to the previous proof and
is omitted. \dokaend

This realization permits one to define a product
\[
K^0\left( \overline{T^{*}M}^\pi \right) \times K_0\left( \mathcal{A}_{M,\pi
}\right) \longrightarrow K_0\left( \mathcal{A}_{T^{*}M,\pi }\right)
\]
using the construction of a symbol with coefficients in a vector bundle. More precisely, for
elements
\[
\left[ \sigma \right] \in K\left( \overline{T^{*}M}^\pi \right) ,\qquad \left[
E,F,\sigma ^{\prime }\right] \in K_0\left( \mathcal{A}_{M,\pi }\right),
\]
consider the symbol
\begin{equation}
\sigma \otimes 1_E\oplus \sigma ^{-1}\otimes 1_F  \label{st2}
\end{equation}
on $M.$ The direct image of its restriction to the boundary can be written as
\[
\begin{array}{ccl}
\pi _{!}\left( \left. \sigma \otimes 1_E\oplus \sigma ^{-1}\otimes
1_F\right| _{\partial M}\right) & = & \pi _{!}\sigma \otimes 1_{\pi
_{!}\left. E\right| _{\partial M}}\oplus \pi _{!}\sigma ^{-1}\otimes 1_{\pi
_{!}\left. F\right| _{\partial M}} \vspace{3mm} \\
&  & \simeq \left( \pi _{!}\sigma \oplus \pi _{!}\sigma ^{-1}\right) \otimes
1_{\pi _{!}\left. E\right| _{\partial M}}.
\end{array}
\]
The latter isomorphism is induced by a vector bundle isomorphism
\[
\pi _{!}\left. E\right| _{\partial M}\stackrel{\sigma ^{\prime }}{\simeq }%
\pi _{!}\left. F\right| _{\partial M}.
\]
Now $\left( \pi _{!}\sigma \oplus \pi _{!}\sigma ^{-1}\right) \otimes 1_{\pi
_{!}\left. E\right| _{\partial M}}$ is trivially homotopic to the identity. The
homotopy is
\[
\left(
\begin{array}{cc}
\pi _{!}\sigma & 0 \\
0 & 1
\end{array}
\right) \left(
\begin{array}{cc}
\cos \tau & \sin \tau \\ %
-\sin \tau & \cos \tau%
\end{array}
\right) \left(
\begin{array}{cc}
1 & 0 \\
0 & \pi _{!}\sigma ^{-1}
\end{array}
\right) \left(
\begin{array}{cc}
\cos \tau & -\sin \tau \\ %
\sin \tau & \cos \tau%
\end{array}
\right) ,\quad \tau \in \left[ 0,\pi /2\right] .
\]
Thus we have extended the symbol (\ref{st2}) to a non-local elliptic symbol on
$M.$ Now the desired product is defined as the difference construction of the
latter symbol, which we denote by
\[
\left[ \sigma \right] \times \left[ E,F,\sigma ^{\prime }\right] \in
K_0\left( \mathcal{A}_{T^{*}M,\pi }\right) .
\]
Using the homotopy classification  $K_0\left( \mathcal{A}_{T^{*}M,\pi
}\right)\simeq {\rm Ell}(M,\pi)$, we can apply the index mapping to this
product and define the pairing of groups as the composition
\[
\left\langle ,\right\rangle :K^0\left( \overline{T^{*}M}^\pi \right) \times
K_0\left( \mathcal{A}_{M,\pi }\right) \longrightarrow K_0\left( \mathcal{A}%
_{T^{*}M,\pi }\right) \stackrel{\limfunc{ind}}{\longrightarrow }\Bbb{Z}.
\]
This pairing is an analogue of Poincar\'e duality on the singular manifold
 $\overline{T^{*}M}^\pi $
(see Section~\ref{dvoj6}). Let us note that for a regular covering the index
can be computed topologically (by the index theorem) and hence the pairing is
also topologically computable. \vspace{1mm}

\textbf{3}. \textbf{The element of }$K$\textbf{-theory with coefficients
defined by a manifold with a covering on the boundary.} We start from a
universal example. Denote the half-infinite cylinder $EG_N\times [0,+\infty)$
by  ${\cal M}_N$, and the projective limit of the groups
$K_0({\cal{A}}_{{\cal{M}}_N,\pi_N})$ as $N\to \infty$ by
$K_0({\cal{A}}_{{\cal{M}}_\infty})$. Let us show that the universal bundle
$\gamma=\pi_!1 \in \limfunc{Vect}\left( BG\right), $ where $1\in{\rm
Vect}(EG)$, defines an element
\begin{equation}\label{volna}
[\widetilde{\gamma}]\in K_0(\mathcal{A}_{\mathcal{M}_\infty},\mathbb{Q}/n\mathbb{Z}),\quad n=|G|.
\end{equation}
To this end, note that the universal bundle gives the element\footnote{ Here
and below, the $K$-groups of classifying spaces are defined as the projective
limits over their finite-dimensional approximations.}
\[
\left[ \gamma \right] -n\in K^0(BG).
\]
We use the following lemma to show that this difference defines the desired element.
\begin{lemma}
\label{lu2}  There is an isomorphism $\widetilde{K}^0\left( BG\right) \simeq K^1\left( BG,\Bbb{Q}/n\Bbb{Z}%
\right),$ which is defined as the coboundary mapping  $\partial$ in the exact
sequence
\[
\rightarrow K^1\left( BG _N\right) \otimes \Bbb{Q}%
\longrightarrow K^1\left(  BG _N,\Bbb{Q}/n\Bbb{Z}\right) \stackrel{\partial}\longrightarrow
K^0\left( BG _N\right)
\stackrel{\times n}{%
\longrightarrow }K^0\left( BG _N\right) \otimes \Bbb{Q}
\]
induced by the inclusion of the coefficient groups $n\Bbb{Z\subset Q}$.
\end{lemma}

\noindent\emph{Proof}. Let us rewrite the sequence
\[
\rightarrow K^1\left(  BG _N\right) \otimes \Bbb{Q}%
\longrightarrow K^1\left( BG_N,\Bbb{Q}/n\Bbb{Z}\right)
\longrightarrow K^0\left( BG_N\right) \stackrel{\times n}{%
\longrightarrow }K^0\left( BG _N\right) \otimes \Bbb{Q}
\]
as the short exact sequence
\[
0\rightarrow K^1\left(  BG _N\right) \otimes \Bbb{Q}/n\Bbb{Z}%
\longrightarrow K^1\left(  BG _N,\Bbb{Q}/n\Bbb{Z}\right) \longrightarrow \limfunc{Tor}K^0\left(
BG _N\right) \rightarrow 0.
\]
Then  we obtain the following  sequence of projective limits as $N\rightarrow
\infty $ (this sequence may not be exact):
\begin{equation}
0\rightarrow \stackunder{\longleftarrow }{\lim }K^1\left(  BG
_N\right) \otimes \Bbb{Q}/n\Bbb{Z}\longrightarrow K^1\left( BG,\Bbb{Q}/n%
\Bbb{Z}\right) \longrightarrow \,\stackunder{\longleftarrow }{\lim }\limfunc{%
Tor}K^0\left(  BG _N\right) \rightarrow 0.  \label{injlim1}
\end{equation}
As $N$ increases, the sequence $\widetilde{K}^*\left(  BG _N\right) $ has the
following property (e.g., see \cite{APS2,Ati6}). For an arbitrary $N$,  there
exists an $L>0$ such that the range of the mapping
\begin{equation}
\widetilde{K}^*\left(  BG _{N+L}\right) \longrightarrow \widetilde{K}^*\left( BG _{N}\right)
\label{svo1}
\end{equation}
is in the torsion subgroup.  Using this, we obtain the following expressions
for the limits:
\[
\stackunder{\longleftarrow }{\lim }K^1\left(  BG _N\right)
\otimes \Bbb{Q}/n\Bbb{Z}=0,\quad\stackunder{\longleftarrow }{\lim }%
\limfunc{Tor}K^0\left(  BG _N\right) =\,\stackunder{%
\longleftarrow }{\lim }\widetilde{K}^0\left( BG _N\right) =%
\widetilde{K}^0\left( BG\right) .
\]
One can also use (\ref{svo1}) to prove the exactness of  (\ref{injlim1}). The
proof is based on commutative diagrams of the form
\begin{equation}
\begin{array}{ccccc}
0\rightarrow K^1\left(  BG _N\right) \otimes \Bbb{Q}/n\Bbb{Z}
& \!\!\longrightarrow \!\! & K^1\left(  BG _N,\Bbb{Q}/n\Bbb{Z}%
\right) & \!\!\longrightarrow \!\! & \limfunc{Tor}K^0\left(
BG _N\right) \rightarrow 0 \\
0\uparrow &  & \uparrow &  & \uparrow \\
0\rightarrow \!K^1\!\left(  BG _{N+L}\right) \otimes \Bbb{Q}/n%
\Bbb{Z} & \!\!\longrightarrow \!\! & \!K^1\!\left( BG
_{N+L}\!,\Bbb{Q}/n\Bbb{Z}\right) \! & \!\!\longrightarrow \!\! & \!\limfunc{%
Tor}\!K^0\!\left(  BG _{N+L}\right) \!\rightarrow 0,
\end{array}
\label{coma1}
\end{equation}
where  $L$ is chosen as in  (\ref{svo1}).

Thus we obtain the desired isomorphism
\[
\stackunder{\longleftarrow }{\lim }K^1\left(  BG _N,\Bbb{Q}/n%
\Bbb{Z}\right) \simeq \stackunder{\longleftarrow }{\lim }\widetilde{K}%
^0\left( BG _N\right) .
\]
\dokaend

Finally, the desired element
 \eqref{volna} is obtained from the isomorphism
$$
K^{*+1}(BG_N,\Bbb{Q}/n\Bbb{Z})\simeq K_*({\cal{A}}_{{\cal{M}}_N,\pi_N}, \Bbb{Q}/n\Bbb{Z})
$$
induced by the inclusion of the ideal
 $C_0(BG_N\times (0,1),{\rm End}(\pi_N)_!1)$. Thus
$$ K^1\left( BG,\Bbb{Q}/n\Bbb{Z}\right) \simeq
K_0({\cal{A}}_{{\cal{M}}_\infty},\Bbb{Q}/n\Bbb{Z}),
$$
and $\left[ \widetilde{\gamma }\right]$ can be viewed as an element of both
groups.

Suppose that we are now given a pair $\left( M,\pi \right).$ In the remaining
part of the section, we assume that $\pi$ is a principal  $G$-covering for a
finite group $G$. There exists a mapping $f:M\to{\cal M}_N$ that takes the
boundary to the base $EG_N\times\{0\}$ and is equivariant on the boundary (see
Proposition~\ref{pre1}). The inverse image of $[\widetilde{\gamma}]$ is denoted
by
\begin{equation}
[\widetilde{\pi_!1}]\stackrel{\rm def}={f}^{*}[ \widetilde{\gamma }] \in K_0\left( \mathcal{A}%
_{M,\pi },\Bbb{Q}/n\Bbb{Z}\right).  \label{navor1}
\end{equation}
This element does not depend on the choice of $f$, since a map into the
universal space is unique up to homotopy. Let us obtain a geometric realization
of this element. We do this in two steps.

\textbf{4. A geometric realization of $K_0\left( \mathcal{A}%
_{M,\pi },\Bbb{Q}/n\Bbb{Z}\right)$.} The group $\Bbb{Q}/n\Bbb{Z}$ is the direct
limit of the finite groups
\[
\Bbb{Z}_{nN}  \subset  \Bbb{Q}/n\Bbb{Z}, \quad x \mapsto  x/N.
\]
Therefore, the $K$-group with coefficients in $\Bbb{Q}/n\Bbb{Z}$ is defined as
the direct limit
\begin{equation}\label{alp1}
K_0\left( \mathcal{A}_{M,\pi },\Bbb{Q}/n\Bbb{Z}\right) =\stackunder{%
\longrightarrow }{\lim }K_0\left( \mathcal{A}_{M,\pi },\Bbb{Z}_{nN}\right).
\end{equation}

Further, the elements of the groups with finite coefficients can be constructed
by using the following proposition (cf. \cite{SaScS3}).
\begin{proposition}\label{pre5}
A triple $\left( E,F,\sigma \right) ,$ where
\begin{equation}
E\in \limfunc{Vect}\left( M\right) ,\;F\in \limfunc{Vect}\left( X\right) ,\qquad \pi _{!}\left(
\left. E\right| _{\partial M}\right) \stackrel{\sigma }{\simeq }kF,  \label{modnn1}
\end{equation}
and $\sigma $ is an isomorphism on $X$, defines an element in $K_0\left(
\mathcal{A}_{M,\pi },\Bbb{Z}_k\right) .$
\end{proposition}

\noindent\emph{Proof.} By analogy with the topological case (e.g., see
\cite{APS2}), the theory with coefficients in  $\Bbb{Z}_k$ is defined in terms
of the \emph{Moore space} $\mathbb{M}_k$ of this group by the formula
\begin{equation}
K_0\left( \mathcal{A}_{M,\pi },\Bbb{Z}_k\right) =K_0\left( \widetilde{C}%
_0\left( \mathbb{M}_k,\mathcal{A}_{M,\pi }\right) \right) , \label{nesht1}
\end{equation}
where $\widetilde{C}_0\left( \mathbb{M}_k,\mathcal{A}_{M,\pi }\right) $ is the algebra of
$\mathcal{A}_{M,\pi }$-valued functions on the Moore space vanishing at the fixed point.

One can readily generalize Lemma~\ref{lemc2} to the case of families. More
precisely, the same method shows that the group $K_0\left(
\widetilde{C}_0\left( \mathbb{M}_k,\mathcal{A}_{M,\pi }\right) \right) $ is
isomorphic to the group of stable homotopy classes of triples $\left( E^{\prime
},F^{\prime },\sigma ^{\prime }\right) ,$ where $E^{\prime },F^{\prime }\in
\limfunc{Vect}\left( M\times \mathbb{M}_k\right) $ and the isomorphism
\[
\sigma ^{\prime }:\pi _{!}\left(\left. E^{\prime }\right| _{\partial M}\right)\longrightarrow \pi
_{!}\left(\left. F^{\prime }\right| _{\partial M}\right)
\]
is defined over $X\times \mathbb{M}_k.$

Let $\varepsilon$ be the line bundle over the Moore space representing the
generator
$\left[\varepsilon\right]-1\in\widetilde{K}\left(\mathbb{M}_k\right)\simeq\Bbb{Z}_k$.
(Further information about the Moore spaces can be found, e.g., in
\cite{APS2,SaScS3}.) Let us also fix a trivialization $\rho :k\varepsilon
\rightarrow \Bbb{C}^k.$

To the triple $\left( E,F,\sigma \right)$ in \eqref{modnn1}, we assign  the
element
\[
\left[ E\otimes \varepsilon ,E,\sigma ^{\prime }\right] \in K_0\left( \widetilde{C}_0\left(
\mathbb{M}_k,\mathcal{A}_{M,\pi }\right) \right) ,
\]
where the isomorphism  $\sigma ^{\prime }$ is defined as the composition (see \cite{SaScS3})
\begin{equation}
\pi _{!}\left(\left. E\right| _{\partial M}\right)\otimes \varepsilon \stackrel{\sigma \otimes
1}{\rightarrow }kF\otimes \varepsilon \simeq F\otimes k\varepsilon
\stackrel{1\otimes \rho }{\rightarrow }F\otimes \Bbb{C}^k\simeq kF\stackrel{%
\sigma ^{-1}\otimes 1}{\rightarrow }\pi _{!}\left(\left. E\right| _{\partial M}\right).
\label{longo1}
\end{equation}
\dokaend

\textbf{5.  A geometric  realization of $[\widetilde{\pi_! 1}]$
(see~\eqref{navor1}).} For sufficiently large $N$, consider a triple $\left(
N,1,\alpha \right) ,$ where
\[
N\in \limfunc{Vect}\left( M\right) ,\qquad 1\in \limfunc{Vect}\left( X\right) ,
\]
are trivial vector bundles of the corresponding dimensions and $ \pi
_{!}N\stackrel{\alpha }{\simeq }\Bbb{C}^{nN}$ is some trivialization. By
Proposition~\ref{pre5}, this triple defines an element
\[
\left[ N,1,\alpha \right] \in K_0\left( \mathcal{A}_{M,\pi },\Bbb{Q}/n\Bbb{Z}%
\right).
\]
Now, using the diagram (\ref{coma1}), the reader can verify that this element
coincides with $[ \widetilde{\pi_{!}1}]$ if the number $N$ and the
trivialization are chosen as follows.

Suppose that the range of the classifying mapping
\[
f:M\longrightarrow \mathcal{M}_\infty
\]
is contained in the skeleton $\mathcal{M}_{N^{\prime }},$ then for $N^{\prime
}$ there exists an $L^{\prime }$ such that property~(\ref{svo1}) is valid. Now
we can choose an $N$ such that the restriction of the direct sum $N\gamma $ of
the universal bundle  to $BG _{N^{\prime }+L^{\prime }}$ is trivial with some
trivialization
\[
N\gamma \stackrel{\alpha ^{\prime }}{\simeq }\Bbb{C}^{nN}.
\]
Finally, over $M$ we choose the induced trivialization
$$
\alpha=f^*\alpha'.
$$

\textbf{6. The index defect theorem}.

\begin{theorem}
\label{defect2}Let $(M,\pi)$  be a manifold with a covering on the boundary
corresponding to a free action of a finite group $G$. Then the diagram
$$
\xymatrix{\limfunc{Ell}\left( \overline{M}^\pi \right) \ar[d]_{\widetilde{\limfunc{ind}}}
\ar[r]^\chi & K\left( \overline{T^{*}M}^\pi \right) \ar[ld]^{\left\langle
\cdot , [\widetilde{\pi _{!}1}] \right\rangle} \\
\Bbb{R}/n\Bbb{Z}\text{,} }
$$
commutes. Here $\left\langle ,\right\rangle $ is the Poincar\'e pairing with
coefficients,
\begin{equation}
\left\langle ,\right\rangle :K\left( \overline{T^{*}M}^\pi \right) \times K_0\left(
\mathcal{A}_{M,\pi },\Bbb{Q}/n\Bbb{Z}\right) \longrightarrow
K_0\left( \mathcal{A}_{T^{*}M,\pi },\Bbb{Q}/n\Bbb{Z}\right) \stackrel{%
\limfunc{ind}}{\rightarrow }\Bbb{Q}/n\Bbb{Z}.  \label{proda2}
\end{equation}
\end{theorem}

\begin{remark}
\emph{Theorem \ref{defect2} expresses $\widetilde{\rm ind}D$ in topological terms via the
principal symbol. Indeed, by  \eqref{alp1} and \eqref{nesht1}, the index mapping in
\eqref{proda2} can be expressed topologically using the index theorem for families (Theorem
 \ref{indfam}).}
\end{remark}

\noindent\emph{Proof}. The proof of the theorem is essentially analytic in nature. The main idea
is to reduce the analytic invariant $\widetilde{\rm ind}$ to the index with values in
$\mathbb{Q}/n\mathbb{Z}$. Thus, we start by defining the corresponding operators.

1. First, we define elliptic theory  $\limfunc{Ell}\left( M,\pi ,\Bbb{Q}/n\Bbb{Z%
}\right) $ with coefficients $\Bbb{Q}/n\Bbb{Z}$. The definition can be given using the direct
limit
\[
\limfunc{Ell}\left( M,\pi ,\Bbb{Q}/n\Bbb{Z}\right) =\stackunder{%
\longrightarrow }{\lim }\limfunc{Ell}\left( M,\pi ,\Bbb{Z}_{nN}\right) ,\quad \Bbb{Z}_{nN}\subset
\Bbb{Z}_{nNM}\subset \Bbb{Q}/n\Bbb{Z},
\]
of theories with finite coefficients. More precisely, elliptic theory with
coefficients in $\Bbb{Z}_k$ is defined by families of non-local elliptic
operators of order one parametrized by the Moore space $\mathbb{M}_k$ of
$\Bbb{Z}_k$:
\[
\limfunc{Ell}\left( M,\pi ,\Bbb{Z}_k\right) =\limfunc{Ell}%
_{\mathbb{M}_k}\left( M,\pi \right).
\]
For  elliptic theory with coefficients, we refer the reader to  \cite{SaScS3}.

2. Consider the mapping
\begin{equation}\label{phibig}
\begin{array}{ccc}
\limfunc{Ell}\left( \overline{M}^\pi \right) & \stackrel{\Phi }{%
\longrightarrow } & \limfunc{Ell}\left( M,\pi ,\Bbb{Q}/n\Bbb{Z}\right)
\end{array}
\end{equation}
that takes an operator $D$ to the  family
\[
D^{*}\oplus \left( D\otimes 1_\varepsilon \right) :C^\infty \left( M,F\oplus E\otimes \varepsilon
\right) \longrightarrow C^\infty \left( M,E\oplus F\otimes \varepsilon \right)
\]
of first-order elliptic operators on $M$ parametrized by $\mathbb{M}_{nN}$ (the
number $N$ will be chosen below). Here $D^{*}$ is the adjoint operator, and the
operator family obtained by twisting $D$ with the bundle $\varepsilon $ is
denoted by $D\otimes 1_\varepsilon $. Consider the direct sum of $N$ copies of
this family. It turns out that if $N$ is sufficiently large, then this family
admits an elliptic boundary condition. Indeed, for $N$ sufficiently large there
exists a trivialization
\begin{equation}
N\pi _{!}1\stackrel{\alpha }{\simeq }\Bbb{C}^{nN},  \label{triva1}
\end{equation}
since $\pi_!1$ is flat. Hence on the base of the covering we have the vector
bundle isomorphism
\[
\pi _{!}\left( N\left. E\right| _{\partial M}\right) \simeq \pi _{!}N\otimes E_0\stackrel{\alpha
\otimes 1}{\longrightarrow }\Bbb{C}^{nN}\otimes E_0
\]
and the similar isomorphism
\[
\pi _{!}\left( N\left. E\right| _{\partial M}\right) \otimes \varepsilon
\simeq \pi _{!}N\otimes E_0\otimes \varepsilon \stackrel{\alpha \otimes 1}{%
\longrightarrow }\Bbb{C}^{nN}\otimes E_0\otimes \varepsilon \simeq
nN\varepsilon \otimes E_0\stackrel{\rho \otimes 1}{\longrightarrow }\Bbb{C}%
^{nN}\otimes E_0.
\]
We denote the induced isomorphisms on sections by
\[
B_1=\alpha \otimes 1:C^\infty \left( X,\pi _{!} \left(N\left. E\right| _{\partial
M}\right)\right) \stackrel{}{\longrightarrow }C^\infty \left( X,\Bbb{C}^{nN}\otimes E_0\right) ,
\]
\[
B_2=\rho \otimes 1\left( \alpha \otimes 1\right) :C^\infty \left( X,\pi _{!}
\left(N\left. E\right| _{\partial M}\right)\otimes \varepsilon\right) \stackrel{}{%
\longrightarrow }C^\infty \left( X,\Bbb{C}^{nN}\otimes E_0\right) .
\]
Now we can define the family of non-local boundary value problems (in the sense
of Section~\ref{para1})
\begin{equation}
\left\{
\begin{array}{l}
\begin{array}{cc}
ND^{*}u=f_1,\qquad & N\left( D\otimes 1_\varepsilon \right) v=f_2,
\end{array}
\vspace{2mm} \\
B_1\beta _E\left. u\right| _{\partial M}+B_2\beta _E\left. v\right| _{\partial M}=g,\quad g\in
C^\infty \left( X,\Bbb{C}^{nN}\otimes E_0\right) ,
\end{array}
\right.  \label{nelin1}
\end{equation}
which consists of elliptic elements. We define the mapping \eqref{phibig} as
follows: it takes $D$ to the family of non-local problems \eqref{nelin1}. Note
that $\Phi$ depends on the choice of the trivialization \eqref{triva1}.

3. There is a natural index mapping
\[
\limfunc{ind}:\limfunc{Ell}\left( M,\pi ,\Bbb{Q}/n\Bbb{Z}\right) \longrightarrow \Bbb{Q}/n\Bbb{Z}
\]
that takes an element $\left[ \mathcal{D}\right] $  represented by a family  $\mathcal{D}$ of
elliptic operators parametrized by $\mathbb{M}_{nN}$ to the (reduced) index of the family
\[
\limfunc{ind}\left[ \mathcal{D}\right] \stackrel{\rm def}=\limfunc{ind}\mathcal{D}\in
\widetilde{K}\left( \mathbb{M}_{nN}\right) \simeq \Bbb{Z}_{nN}\subset \Bbb{Q}/n\Bbb{Z}.
\]

\begin{lemma}
\label{lem3} The diagram
$$
\xymatrix{ \qquad\limfunc{Ell}\left( \overline{M}^\pi \right) \ar[rr]^{\Phi}
\ar[rd]_{\widetilde{\limfunc{ind}}} & & \limfunc{Ell}\left( M,\pi ,\Bbb{Q}/n\Bbb{Z}\right)
\ar[ld]^{\rm ind} \\
& \quad \Bbb{R}/n\Bbb{Z}, }
$$
where $N$ and the trivialization $\alpha$ in \eqref{triva1} are chosen as in
Subsec.~\emph{5}, commutes.
\end{lemma}

\noindent\emph{Proof of the lemma}. The boundary value problem (\ref{nelin1})
is linearly homotopic to the problem
\[
\left\{
\begin{array}{l}
ND^{*}u=f_1,\qquad N\left( D\otimes 1_\varepsilon \right) v=f_2,\vspace{2mm}
\\
B_1\beta _E\Pi _{-}\left. u\right| _{\partial M}+B_2\beta _E\Pi _{+}\left.
v\right| _{\partial M}=g,\quad g\in C^\infty \left( X,\Bbb{C}^{nN}\otimes
E_0\right),
\end{array}
\right.
\]
within the class of elliptic problems. The last formula shows that the index of $\Phi \left[
D\right] $ is equal to the sum of the index of the family of spectral problems for $ND^{*}$ and
$N\left( D\otimes 1_\varepsilon \right) $ and the index of the operator family
\begin{equation}
N\func{Im}\Pi _{-}\left( A\right) \oplus N\func{Im}\Pi _{+}\left( A\right)
\otimes \varepsilon \stackrel{B_1+B_2}{\longrightarrow }C^\infty \left( X,%
\Bbb{C}^{nN}\otimes E_0\right)  \label{opcl1}
\end{equation}
on the boundary. Note that we specify the self-adjoint operators in the notation of spectral
projections. Let us compute the index of the former family on $X$.

1) There is a decomposition
\[
C^\infty \left( X,\Bbb{C}^{nN}\otimes E_0\right) \simeq nN\func{Im}\Pi _{-}\left( A_0\right)
\oplus nN\varepsilon \otimes \func{Im}\Pi _{+}\left( A_0\right)
\]
of the  target space for the family (\ref{opcl1}). This decomposition is defined as
\[
nN\func{Im}\Pi _{-}\left( A_0\right) \oplus nN\varepsilon \otimes \func{Im}%
\Pi _{+}\left( A_0\right) \stackrel{1+\left( \rho \otimes 1\right) }{%
\longrightarrow }C^\infty \left( X,\Bbb{C}^{nN}\otimes E_0\right) .
\]
Using this isomorphism, we represent the index of (\ref{opcl1}) in the form
\[
=\limfunc{ind}\left( N\func{Im}\Pi _{+}\left( A\right) \stackrel{\Pi _{+}\left( A_0\right) \beta
_E}{\longrightarrow }nN\func{Im}\Pi _{+}\left( A_0\right) \right) \left( \left[ \varepsilon
\right] -1\right) \in \widetilde{K}\left( \mathbb{M}_{nN}\right) .
\]
Finally, we rewrite the index by pushing forward the space $\func{Im} \Pi
_{+}\left( A\right) $ to the base of the covering:
\[
=\limfunc{ind}\left( N\func{Im}\Pi _{+}\left( \pi _{!}A\right) \stackrel{\Pi _{+}\left(
A_0\right) }{\longrightarrow }nN\func{Im}\Pi _{+}\left( A_0\right) \right) \left( \left[
\varepsilon \right] -1\right) .
\]
The index of the elliptic operator (not a family!) in the last formula can be expressed by the
Atiyah--Patodi--Singer formula \cite{APS3}
\begin{equation}
\limfunc{ind}\left( N\func{Im}\Pi _{+}\left( \pi _{!}A\right) \stackrel{\Pi _{+}\left( A_0\right)
}{\longrightarrow }nN\func{Im}\Pi _{+}\left( A_0\right) \right) =N\eta \left( A\right) -nN\eta
\left( A_0\right) +\left\langle \left[ \sigma \left( A_0\right) \right] ,\left[ \pi _{!}1\right]
\right\rangle , \label{apps}
\end{equation}
where the brackets $\left\langle,\right\rangle$ denote the pairing
\begin{equation}
\left\langle ,\right\rangle :K^1\left( T^{*}X\right) \times K^1\left( X,\Bbb{%
Q}\right) \longrightarrow \Bbb{Q}  \label{pai2}
\end{equation}
of the difference element $ \left[ \sigma \left( A_0\right) \right] \in
K^1\left( T^{*}X\right) $ of an elliptic self-adjoint operator $A_0$ with the
element $\left[ \pi _{!}1\right] \in K^1\left( X,\Bbb{Q}\right) $ defined by
the trivialized flat bundle $N\pi _{!}1$
(more about this formula can be found in the book \cite{Gil1}).%

2) It turns out that for our choice of the trivialization (\ref{triva1}) the
last term in (\ref{apps}) is equal to zero. Indeed, consider the classifying
mapping $f:X\rightarrow BG _{N^{\prime }}.$ We can evaluate
 (\ref{pai2}) on the classifying space:
\begin{equation}
\left\langle \left[ \sigma \left( A_0\right) \right] ,\left[ \pi _{!}1\right] \right\rangle
=\left\langle f_{!}\left[ \sigma \left( A_0\right) \right] ,\left[ \gamma \right] \right\rangle
,\quad \left[ \pi _{!}1\right] =f^{*}\left[ \gamma \right] \in K^1\left( X,\Bbb{Q}\right),
\label{al1}
\end{equation}
where $[\gamma]\in K^1(BG_{N'})\otimes\Bbb{Q}$ is the element defined by the trivialized flat
bundle  $N\gamma$. The inclusion $ BG _{N^{\prime }}\subset BG _{N^{\prime }+L^{\prime }}$
induces the commutative diagram
\[
\begin{array}{ccccc}
K^1\left( T^{*} BG _{N^{\prime }}\right) & \times & K^1\left(  BG _{N^{\prime }}\right) \otimes
\Bbb{Q} & \longrightarrow &
\Bbb{Q} \\
\downarrow &  & \uparrow &  & \parallel \\
K^1\left( T^{*} BG _{N^{\prime }+L^{\prime }}\right) & \times & K^1\left(  BG _{N^{\prime
}+L^{\prime }}\right) \otimes \Bbb{Q} & \longrightarrow & \Bbb{Q}
\end{array}
\]
Using this diagram and (\ref{svo1}), one can prove the triviality of the pairing (\ref{al1}) by a
diagram chase argument.

Thus we have reduced $\limfunc{ind}\Phi \left[ D\right] $ to the desired form
\[
\limfunc{ind}\Phi \left[ D\right] =\widetilde{\limfunc{ind}}\left[ D\right] .
\]
\dokaend

4. To complete the proof of the theorem, it suffices to show that the value of the Poincar\'e
pairing $\langle [\sigma(D)],[\widetilde{\pi_!1}]\rangle$ coincides with the index of $\Phi(D)$.

We denote the product by $[ \widetilde{\pi _{!}1}] $ by $\varphi$:
\[
\varphi :K\left( \overline{T^{*}M}^\pi \right) \longrightarrow K_0\left( \mathcal{A}_{T^{*}M,\pi
},\Bbb{Q}/n\Bbb{Z}\right).
\]
\begin{lemma}
\label{lem3aa}Under the assumptions of Lemma \emph{\ref{lem3}}, the diagram
\begin{equation}
\begin{array}{ccc}
\limfunc{Ell}\left( \overline{M}^\pi \right) & \stackrel{\Phi }{%
\longrightarrow } & \limfunc{Ell}\left( M,\pi ,\Bbb{Q}/n\Bbb{Z}\right) \\
{\chi \downarrow } &  & \downarrow {\chi ^{\prime }} \\ %
K\left( \overline{T^{*}M}^\pi \right) & \stackrel{\varphi }{\longrightarrow } & K\left(
\mathcal{A}_{T^{*}M,\pi },\Bbb{Q}/n\Bbb{Z}\right),
\end{array}
\label{inty3}
\end{equation}
where  $\chi'$ is induced by the difference constructions for families
\emph{(}see Subsec.~{\rm 3.3}{\em )}.
\end{lemma}

\noindent\emph{Proof}. Substituting the definitions of  $\left[ \sigma
\left(D\right) \right] $ and $[\widetilde{\pi _{!}1}] $ (according to
Subsecs.~4.4 and 4.5) into Eq.~(\ref{st2}), defining the product, one can show
that the desired product $\varphi[\sigma(D)]$ is determined by the family of
elliptic symbols that  are equal to
\begin{equation}
N\sigma \left( D\right) \otimes 1_\varepsilon \oplus N\sigma \left( D\right)
^{-1}\otimes 1 \label{aa1}
\end{equation}
far from the boundary. The direct image of the restriction of this symbol to
the boundary is equal to
\begin{eqnarray}
N\pi _{!}\left( \sigma \left( D\right) \otimes 1_\varepsilon \oplus \sigma \left( D\right)
^{-1}\otimes 1\right) &=&\left( \sigma \left( D_0\right) \otimes 1_{N\varepsilon \otimes \pi
_{!}1}\oplus \sigma \left( D_0\right)
^{-1}\otimes 1_{N\pi _{!}1}\right) \simeq  \label{aa2} \\
&\simeq &\left( \sigma \left( D_0\right) \oplus \sigma \left( D_0\right)
^{-1}\right) \otimes 1_{\Bbb{C}^{nN}}. \nonumber
\end{eqnarray}
In the last equality, we use the isomorphisms
$N\pi _{!}1\stackrel{\alpha }{\simeq }\Bbb{C}^{nN},\quad \Bbb{C}%
^{nN}\varepsilon \stackrel{\rho }{\simeq }\Bbb{C}^{nN}.$ The symbol is extended
to a neighbourhood of the boundary  using the homotopy of the direct sum
$\sigma \left( D_0\right)\oplus \sigma \left( D_0\right) ^{-1}$ to the
identity.

It remains to prove that the difference elements for the principal symbol of
the family of boundary value problems (\ref {nelin1}) and the symbol defined by
(\ref{aa1}), (\ref{aa2}) coincide. Indeed, this equality is obvious far from
the boundary, since the only difference here is in the components $\sigma
\left( D^{*}\right) $ and $\sigma \left( D\right) ^{-1}.$ These components are
joined by the standard homotopy
\[
\sigma \left( D^{*}\right) \left[ \sigma \left( D\right) \sigma \left( D^{*}\right) \right]
^{-s},\qquad s\in \left[ 0,1\right] .
\]
The reader can also prove the equality near the boundary using the formulae for
order reduction given in Remark~\ref{remk}. \dokaend

By combining Lemmata~\ref{lem3aa} and \ref{lem3}, we complete the proof of the
theorem. \dokaend

\section{Applications\label{exa4}}

\textbf{1}. Theorem \ref{defect2} enables one to express the fractional part of the $\eta
$-invariant in the following situation.

Let $M$ be an even-dimensional spin manifold with boundary represented as the
total space of a covering such that the spin structure on the boundary is the
pullback of a spin structure on the base. Let us also fix an $E\in
\limfunc{Vect}\left( M\right) $ that  is also pulled back from the base near
the boundary: $\left. E\right| _{\partial M}\simeq\pi^* E_0$. We choose a
metric on $M$ that is a product metric induced by a metric on the base near the
boundary. Finally, we choose a similar connection in $E$.

\begin{proposition}
The Dirac operator $D_M$ on $M$ with coefficients in $E$ satisfies the
assumptions of Theorem \emph{\ref{defect2}}, and the fractional part of the
$\eta $-invariant is equal to
\[
\left\{ \eta \left( D_X\right) \right\} =\frac 1n\left( \int\limits_M%
\widehat{A}\left( M\right) \limfunc{ch}E-\left\langle \left[ \sigma \left( D_M\right) \right] ,[
\widetilde{\pi _{!}1}] \right\rangle \right) \in \Bbb{R}/\Bbb{Z},
\]
where $D_X$ is the self-adjoint Dirac operator on $X$ with coefficients in $E_0$.
\end{proposition}

\emph{Proof}. The formula follows from Theorem \ref{defect2} if we decompose
the index of the spectral problem using the Atiyah--Patodi--Singer formula
\[
\limfunc{ind}\left( D_M,\Pi _{+}\right) =\int\limits_M\widehat{A}\left(
M\right) \limfunc{ch}E-\eta \left( D_{\partial M}\right) .
\]
\dokaend

\vspace{2mm}

\textbf{2}. The invariant $\widetilde{\limfunc{ind}}$ can be effectively computed via Lefschetz
theory. Suppose that $\pi $ is regular, i.e., the boundary is a principal $G$-bundle for a finite
group $G$. Let $D$ be a $G$-invariant elliptic differential operator of order one on $M.$ For
$g\in G$, let $L\left( D,g\right) \in \Bbb{C}$ be the usual contribution to the Lefschetz formula
(see \cite{Don3}) of the fixed point set of the diffeomorphism $g:M\rightarrow M.$

\begin{proposition}
One has
\begin{equation}
\label{lefa1}\widetilde{\limfunc{ind}}D\equiv -\sum_{g\neq e}L\left( D,g\right) \text{ }%
\left( \func{mod}n\right) .
\end{equation}
\end{proposition}

\emph{Proof}. Consider the equivariant index $\limfunc{ind}_g\left( D,\Pi _{+}\right) $ of the
Atiyah--Patodi--Singer problem and the equivariant $\eta $-function (see \cite{Don3}) of the
tangential operator $A$ on the boundary.

Denote by $\left( D,\Pi _{+}\right) ^G$ and $A^G$ the restrictions of the
corresponding operators to the subspaces of $G$-invariant sections. Clearly,
$A^G$ is isomorphic to $A_0$ on $X$. On the other hand, one can express the
usual invariants in terms of their equivariant counterparts:
\[
\limfunc{ind}\left( D,\Pi _{+}\right) ^G=\frac 1{\left| G\right| }\sum_{g\in
G}\limfunc{ind}_g\left( D,\Pi _{+}\right) \qquad \eta \left( A^G\right) =\frac 1{\left| G\right|
}\sum_{g\in G}\eta \left( A,g\right).
\]
These expression follow from elementary character theory. Using them, we write
\[
\widetilde{\limfunc{ind}}D=\limfunc{ind}_e\left( D,\Pi _{+}\right)
-\sum_{g\neq e}\eta \left( A,g\right) .
\]
Let us substitute the expression for the $\eta $-invariant given by the
equivariant Atiyah--Patodi--Singer formula
(see \cite{Don3}) %
\[
-\eta \left( A,g\right) =\limfunc{ind}_g\left( D,\Pi _{+}\right) -L\left(
D,g\right)
\]
into this formula. This gives the desired congruence \eqref{lefa1}:
\[
\widetilde{\limfunc{ind}}D=\left| G\right| \limfunc{ind}\left( D,\Pi
_{+}\right) ^G-\sum_{g\neq e}L\left( D,g\right) .
\]
\dokaend

\section{Poincar\'{e} isomorphisms}

{\bf 1. A closed smooth manifold.} It is well known (see \cite{Ati4,Kas3}  or
the monograph \cite{HiRo1}) that elliptic operators of order zero on a compact
closed manifold define elements in $K$-theory:
\[
\left[ \sigma \left( D\right) \right] \in K^{*}\left( T^{*}M \right) ,\quad \left[ D\right] \in
K^{*}\left( C\left( M\right) \right) \equiv K_{*}\left( M\right).
\]
The latter group is the analytic $K$-homology group, and the grading is odd for
self-adjoint operators and even otherwise. The first element is the difference
element of the operator. To define the second element, we recall that an
elliptic operator $D$ of order zero is a Fredholm operator
\[
D:L^2\left( M,E\right) \longrightarrow L^2\left( M,F\right) ,
\]
where both $L^2$-spaces are modules over $C\left( M\right) $ (the module structure is given by
the pointwise product of functions). In addition, $D$ commutes with the module structure up to
compact operators. Thus, for a self-adjoint  $D$ (of course in this case the bundles coincide)
the pair $\left( L^2\left( M,E\right) ,D\right) $ is an element
\[
\left[ D\right] \in K^1\left( C\left( M\right) \right) .
\]
For a nonself-adjoint $D$, we consider a self-adjoint matrix
operator
\begin{equation}\label{toper}
T=\left(
\begin{array}{cc}
0 & D^{*} \\
D & 0
\end{array}
\right)
\end{equation}
in the  naturally $\Bbb{Z}_2$-graded $C\left( M\right) $-module $L^2\left(
M,E\right) \oplus L^2\left( M,F\right) $. The operator $T$ is odd with respect
to the grading. Hence it defines a $K$-theory element, denoted by
\[
\left[ D\right] \in K^0\left( C\left( M\right) \right) .
\]

{\bf 2. Manifold with boundary.} On the other hand, elliptic operators of order
one on a manifold with  non-empty boundary define similar elements
\[
\left[ \sigma \left( D\right) \right] \in K^i\left( T^{*}M\right) ,\quad
\left[ D\right] \in K_i\left( M\backslash \partial M\right) .
\]
The former is the Atiyah--Singer difference element, and the latter is defined
as follows. Consider an embedding
\[
M\subset \widetilde{M}
\]
of $M$ in some closed manifold $\widetilde{M}$ of the same dimension (e.g., the
double $2M$). Let $\widetilde{D}$ be an arbitrary extension of $D$ to
$\widetilde{M}$. On $\widetilde{M}$, we consider the zero-order operator
\[
\widetilde{F}=\left( 1+\widetilde{D}^{*}\widetilde{D}\right) ^{-1/2}%
\widetilde{D}.
\]
We define the restriction of this operator to $M$ as the bounded operator
\begin{equation}
F=i^{*}\widetilde{F}i_{*}:L^2\left( M,E\right) \longrightarrow L^2\left(
M,F\right) ,  \label{relo}
\end{equation}
where $i_{*}:L^2\left( M\right) \rightarrow L^2( \widetilde{M}%
) $ is the extension by zero and $i^{*}:L^2( \widetilde{M}) \rightarrow
L^2\left( M\right) $ is the restriction operator.

For a symmetric $D$, we find that $F$ satisfies
\[
\begin{array}{c}
F-F^{*}\in \mathcal{K},\quad f\left( F^2-1\right) \in \mathcal{K}, \quad \left[
F,f\right] \in \mathcal{K}
\end{array}
\]
for functions $f\in C_0\left(
M\backslash \partial M\right)$ vanishing on the boundary. Here
 $\mathcal{K}$ is the ideal of compact operators. These relations show that
$F$ defines an element of $K^1\left( C_0\left( M\backslash
\partial M\right) \right) $ (see
\cite{SoTr1}).%

If $D$ is nonself-adjoint, then one considers the matrix as in
Eq.~\eqref{toper}. This defines an element in $K^0\left( C_0\left( M\backslash
\partial M\right) \right) .$

The mappings
\[
\begin{array}{ccc}
K^{*}\left( T^{*}\left( M\backslash \partial M\right) \right)  &
\longrightarrow  & K_{*}\left( M\right) ,\vspace{1mm} \\
K^{*}\left( T^{*}M\right)  & \longrightarrow  & K_{*}\left( M\backslash
\partial M\right) ,\vspace{3mm} \\
\left[ \sigma \left( D\right) \right]  & \mapsto  & \left[ D\right] ,
\end{array}
\]
which take symbols to operators, define \emph{Poincar\`e isomorphisms} on  a smooth manifold $M$
with boundary (e.g., see \cite{Kas3}). The top mapping is defined in terms of elliptic operators
of order zero on $M$ that are induced by vector bundle isomorphisms near the
boundary.\vspace{1mm}

{\bf 3. Manifolds with singularities.} An elliptic non-local zero-order
operator $D$ defines elements
\[
\left[ \sigma \left( D\right) \right] \in K_{*}\left( \mathcal{A}%
_{T^{*}M,\pi }\right) ,\left[ D\right] \in K^{*}\left( C\left( \overline{M}%
^\pi \right) \right) \simeq K_{*}\left( \overline{M}^\pi \right).
\]
The first is the difference element defined in Section \ref{hom2}. To define
the second element, we note that a non-local elliptic operator $D$ of order
zero does not almost commute with the entire algebra $C\left( M\right)$ but
only with the functions pulled back from the quotient space $\overline{M}^\pi
.$ This leads to a smaller algebra.

On the other hand, the operators of Sections \ref{pr1} and \ref{defe3} define
similar elements
\[
\left[ \sigma \left( D\right) \right] \in K^{*}\left(
\overline{T^{*}M}^\pi \right) ,\left[ D\right] \in K_{*}\left(
\mathcal{A}_{M,\pi }\right) .
\]
In this case, the corresponding operators (\ref{relo}), on the contrary, almost
commute with functions $C_0\left( M\backslash
\partial M\right)$ as well as with the elements of the algebra
$\mathcal{A}_{M,\pi }.$

\begin{theorem}
For an arbitrary manifold with a covering on the boundary $(M,\pi)$, the
following Poincar\'{e} isomorphisms are valid:
\[
\begin{array}{ccc}
K_{*}\left( \mathcal{A}_{T^{*}M,\pi }\right)  & \longrightarrow  &
K_{*}\left( \overline{M}^\pi \right) ,\vspace{1mm} \\
K^{*}\left( \overline{T^{*}M}^\pi \right)  & \longrightarrow  & K^{*}\left(
\mathcal{A}_{M,\pi }\right) ,\vspace{3mm} \\
\left[ \sigma \left( D\right) \right]  & \mapsto  & \left[ D\right] .
\end{array}
\]
\end{theorem}

\emph{Proof}. 1) Let us prove the latter isomorphism. Consider the ideal

\[
I=C_0\left( T^{*}\left( X\times \left( 0,1\right) \right) ,\limfunc{End}%
p^{*}\pi _{!}1\right)\subset \mathcal{A}_{T^{*}M,\pi }
\]
with the quotient $\mathcal{A}_{T^{*}M,\pi }/I\simeq C_0\left( T^{*}M\right).$ The long exact
sequence of the pair can be written as
\begin{equation}
\rightarrow K\left( T^{*}X\right) \stackrel{\alpha }{\rightarrow }K_0\left(
\mathcal{A}_{T^{*}M,\pi }\right) \rightarrow K\left( T^{*}M\right) \rightarrow K^1\left(
T^{*}X\right) \rightarrow \ldots  \label{exa}
\end{equation}
Here we have taken into account the isomorphism $K_{*}\left( C_0\left(
Y,\limfunc{End}G\right) \right) \simeq K_{*}\left( C_0\left( Y\right) \right)
\simeq K^{*}\left( Y\right)$ for a vector bundle
$G\in\limfunc{Vect}\left(Y\right)$.

Consider the commutative diagram
\[
\begin{array}{ccccccccc}
\rightarrow & K^0\left( T^{*}X\right) & \rightarrow & K_0\left( \mathcal{A}%
_{T^{*}M,\pi }\right) & \rightarrow & K^0\left( T^{*}M\right) & \rightarrow
& K^1\left( T^{*}X\right) & \ldots \\
& \downarrow &  & \downarrow &  & \downarrow &  & \downarrow &  \\

\rightarrow & K_0\left( X\right) & \rightarrow & K_0\left( \overline{M}^\pi
\right) & \rightarrow & K_0\left( M,\partial M\right) & \rightarrow &
K_1\left( X\right) & \ldots
\end{array}
\]
Here the lower sequence is the exact sequence of the pair $X\subset
\overline{M}^\pi $ in $K$-homology. The vertical mappings of the diagram
(except for the second one) are isomorphisms (see \cite{BaDo1,Kas3}). Thus,
using the 5-lemma, we find that the middle mapping
\[
K_{*}\left( \mathcal{A}_{T^{*}M,\pi }\right) \longrightarrow K_{*}\left(
\overline{M}^\pi \right)
\]
is also an isomorphism.

2)  In the second case, the proof follows the same scheme, but one uses the diagram
\[
\begin{array}{ccccccccc}
\!\!\leftarrow\!\! & K^1\left( T^{*}X\right) & \!\!\leftarrow\!\! & K^0\left( \overline{%
T^{*}M}^\pi \right) & \!\!\leftarrow\!\! & K^0\left( T^{*}\left( M\backslash
\partial M\right) \right) & \!\!\leftarrow\!\! & K^0\left( T^{*}X\right) & \ldots \\
& \downarrow &  & \downarrow &  & \downarrow &  & \downarrow &  \\
\!\!\leftarrow\!\! & K_1\left( X\right) & \!\!\leftarrow\!\! & K^0\left( \mathcal{A}_{M,\pi
}\right) & \!\!\leftarrow\!\! & K_0\left( M\right) & \!\!\leftarrow\!\! & K_0\left( X\right)
& \ldots
\end{array}
\]
The upper row corresponds to the pair $\Bbb{R}\times T^{*}X\subset \overline{T^{*}M}%
^\pi $.

The proof of the theorem is complete. \dokaend

\section{Poincar\'{e} duality\label{dvoj6}}

An analogue of the pairing for the groups $K^0\left( \overline{T^{*}M}^\pi
\right) $ and $K_0\left( \mathcal{A}_{M,\pi }\right) $ in Section \ref{defe3}
is also valid for the odd groups. The definition is left to the reader.

\begin{theorem}
On a manifold $M$ with covering $\pi$ on the boundary, the pairings
\begin{equation}
K^i\left( \overline{T^{*}M}^\pi \right) \times K_i\left( \mathcal{A}_{M,\pi }\right)
\longrightarrow \Bbb{Z},\qquad i=1,2,  \label{aa}
\end{equation}
are non-degenerate on the free parts of the groups.
\end{theorem}

\emph{Proof}. Fixing the first argument of the pairing, we obtain
a mapping
\[
K^i\left( \overline{T^{*}M}^\pi \right) \otimes \Bbb{Q}\longrightarrow
K_i^{\prime }\left( \mathcal{A}_{M,\pi }\right) ,
\]
where for brevity we write $G^{\prime }=\limfunc{Hom}\left( G,\Bbb{Q}\right) .$
This mapping is part of the commutative diagram
\[
\begin{array}{ccccccccc}
& K^1\left( T^{*}X\right) \otimes \Bbb{Q} & \!\!\leftarrow\!\!  & K^0\left(
\overline{T^{*}M}^\pi \right) \otimes \Bbb{Q} & \!\!\leftarrow\!\!  & K^0\left(
T^{*}\left( M\backslash \partial M\right) \right) \otimes \Bbb{Q} &
\!\!\leftarrow\!\!  & K^0\left( T^{*}X\right) \otimes \Bbb{Q} &  \\
& \downarrow  &  & \downarrow  &  & \downarrow  &  & \downarrow  &  \\
& K^{1\prime }\left( X\right)  & \!\!\leftarrow\!\!  & K_0^{\prime }\left( \mathcal{A%
}_{M,\pi }\right)  & \!\!\leftarrow\!\!  & K^{0\prime }\left(
M\right)  & \!\!\leftarrow\!\! & K^{0\prime }\left( X\right).  &
\end{array}
\]
Here the vertical mappings, except for the second one, are isomorphisms (by
virtue of Poincar\'{e} duality on a closed manifold and on a manifold with
boundary). Thus, by the 5-lemma, the second mapping is an isomorphism. Hence
the pairing (\ref{aa}) is non-degenerate in the second variable.

The non-degeneracy with respect to the first argument can be proved in a
similar way. \dokaend

By way of example, consider $M$ with a $spin^c$-structure that on the boundary
is induced by a $spin^c$-structure on the base $X$ of the covering $\pi .$ Then
the group $K^{*}\left( \overline{T^{*}M}^\pi \right) $ is a free $K^{*+n}\left(
\overline{M}^\pi \right) $-module with one generator (where $n=\dim M$); as a
generator one can take the difference construction
\[
\left[ \sigma \left( D\right) \right] \in K^n\left( \overline{T^{*}M}^\pi
\right)
\]
of the principal symbol of the Dirac operator  on $M$ (this can be proved by
analogy with the usual case of closed manifolds; e.g., see \cite{LaMi1}).
Consequently, one can define the Poincar\'{e} duality pairing
\[
K^{*+n}\left( \overline{M}^\pi \right) \times K_{*}\left( \mathcal{A}_{M,\pi
}\right) \longrightarrow \Bbb{Z}
\]
as the composition with $K^{*+n}\left( \overline{M}^\pi \right) \rightarrow
K\left( \overline{T^{*}M}^\pi \right) .$ The above theorem shows that this
pairing is non-degenerate on the free parts of the groups.

\end{document}